\documentclass[11pt,letterpaper]{amsart}
\usepackage{amssymb,amsfonts,amsthm,array,epsfig,graphics,graphicx,tabularx}
\usepackage{hyperref}
\usepackage[all]{hypcap}
\usepackage{subfigure}
\usepackage{enumerate}

\usepackage[usenames,dvipsnames]{color}

\newenvironment{proo}[1][Proof]{\noindent {\bf #1~: }}{\hfill$\Box$\medskip}

\newtheorem{theorem}{Theorem}[section]
\newtheorem{prop}[theorem]{{Proposition}}

\newtheorem{lemm}[theorem]{{Lemma}}

\newtheorem{clai}[theorem]{{Claim}}
\newtheorem{question}[theorem]{{Question}}

\newtheorem{rema}[theorem]{Remark}

\theoremstyle{definition}
\newtheorem{defi}[theorem]{{Definition}}

\theoremstyle{remark}

\newtheorem{example}[theorem]{{Example}}

\usepackage{fancyhdr}

\title{The indexed links of Non-singular Morse-Smale flows on graph manifolds}
\author{Fangfang Chen and Bin Yu}

\begin{document}
\maketitle

\begin{abstract}
We classify the indexed links 
corresponding to the union of the closed orbits 
of non-singular Morse-Smale flows on  most graph manifolds.
We find that each of this kind of indexed links can be obtained by applying  a finite steps of  operations on a special indexed link, which  consists  of all of the singular Seifert fibers and some regular Seifert fibers with some precisely described conditions.

\end{abstract}

\section{Introduction}\label{s.int}
\subsection{Historic remarks and the aim of the paper}
A \emph{Morse-Smale flow}   is a smooth flow whose chain recurrent set consists of finitely many  hyperbolic closed orbits and fixed points, and it satisfies the transversality condition \cite{Sm}.  If a Morse-Smale flow has no fixed point, then we call it a \emph{non-singular Morse-Smale flow}, abbreviated as an \emph{NMS flow}.

From  the viewpoint of dynamical systems,
Morse-Smale flow is always regarded as a kind of simple system: there does not exist a homoclinic orbit
in a Morse-Smale flow. Roughly speaking, such a system does not provide chaos.
But from the viewpoint more close to topology, i.e.  classifying Morse-Smale flows
up to topological equivalence, it is quite complicated.
Peixoto \cite{Pe} began to systematically classify Morse-Smale flows on surfaces.
 Asimov (\cite{As1}, \cite{As2}) did some significant works from the viewpoint close to topology.
Similar to the well-known relationship between handle decompositions and gradient-like flows
(Morse-Smale flows without closed orbits), he showed that there implies a combinatorial decomposition in an NMS flow, named
a  \emph{round handle decomposition} (abbreviated as \emph{RH decomposition}). Further he used RH decomposition to  obtain several
significant results about NMS flows on $n$-manifolds ($n\geq 4$), one of them (\cite{As1}) says that a closed $n$-manifold admits an NMS flow
if and only if the Euler number of $M$ is zero. This result perfectly answered Question \ref{q.wmad1} in the cases  $n\geq 4$. 
The recent progress about NMS flows we refer to \cite{PS}, \cite{CV3}, \cite{CV}, and \cite{CV2}.

\begin{question}\label{q.wmad1}
Which closed $n$-manifolds admit NMS flows?
\end{question}
Note that by Poincare-Hopf Theorem, it is easy to know that the answer to this question in the case $n=2$ is similar to the cases $n\geq 4$: a closed  surface $\Sigma$ admits an NMS flow if and only if
$\Sigma$ is homeomorphic to either a torus or a Klein bottle. So, for Question \ref{q.wmad1},
what left is the case $n=3$. 
Yano (\cite{refY}) studied the problem of the existence of NMS flows on a closed \(3\)-manifold
in a given homotopy class.
  Morgan (\cite{Mo}) built three theorems to nearly describe the 3-manifolds admitting
NMS flows. One of his theorems says that an irreducible closed orientable $3$-manifold $M$ admits an NMS flow if and only if $M$ is a graph manifold. As a direct consequence, every hyperbolic closed $3$-manifold does not admit any NMS flow.
Note that the Euler number of every closed $3$-manifold is zero, so the case $n=3$  is very special for Question \ref{q.wmad1}.
Therefore it is an interesting topic to more deeply understand NMS flows on $3$-manifolds. In particular, it is natural to ask:
\begin{question}
\label{q.wmad2}
For a given closed orientable $3$-manifold $M$, how to describe the NMS flows on $M$?
\end{question}

The union of closed orbits $\Gamma$ of an NMS flow $\phi_t$ on $M$ is the set of finitely many pairwise disjoint embedded simple closed curves in $M$, i.e.,  a link in $M$. 
In \cite{refK},  Kobayashi described the link type of $\Gamma$ when \(M\) is an irreducible, simple, closed, 
orientable 3-manifold.
We can label every closed orbit $\gamma$ of $\Gamma$
by an integer $k\in \{0,1,2\}$ corresponding to the dimension of the (strong) unstable manifold of $\gamma$, and the labeled link $\Gamma$ is called
the \emph{indexed link} of \(\phi_t\). 
We call this integer $k$ the \emph{index} of $\gamma$.
The indexed link is a  natural dynamical invariant of NMS flows, and is also closely related to knot theory. So, to understand the indexed links of NMS flows on $M$ is a suitable refinement of Question \ref{q.wmad2}.
In a broad sense, 
an \emph{indexed link} in \(M\) is a link in \(M\) with index \(0\), \(1\) or \(2\) attached to each component.
 In \cite{Wa},  Wada  built an algorithm to decide which indexed links can be realized as the union of
closed orbits of NMS flows on $S^3$. This work can be regarded as a solid progress on
Question \ref{q.wmad2} in the case $M\cong S^3$.

The main purpose of this paper is to generalize Wada's description to most graph manifolds. That is to classify the indexed links of NMS flows on  most graph manifolds, which can be regarded as a further progress to answer Question  \ref{q.wmad2}.

 \subsection{Some concepts and notations}
\label{Def}
To state our main results, we have to introduce some further concepts and notations.
\subsubsection{Ordinary graph manifolds}\label{sss.gegra}

We use the notation \(M( \pm g, b; \frac{q_1}{p_1}, \cdots, \frac{q_k}{p_k})\) to record a Seifert fibering of a Seifert  manifold \(M\), where \(g\) is the genus of the base orbifold \(B\), with sign \(+\) if \(B\) is orientable and \(-\) if \(B\) is nonorientable, and \(b\) is the number of the boundary components of \(B\). Here `genus' for nonorientable surfaces means the number of \(\mathbb{RP}^2\) connected summands. Moreover, \(p_i\) and \(q_i\) are coprime, and \(p_i > 0\) for \(i=1, \cdots, k\).
For some \(i_0 \in\{1, \cdots, k\}\), if \(p_{i_0} > 1\), then we say that there is a singular fiber of \(M\) with slope-\(\frac{q_{i_0}}{p_{i_0}}\).  The details we refer to Chapter 2 of Hatcher \cite{Ha}.


A compact irreducible orientable 3-manifold  \(W\) is a \emph{graph manifold} if every JSJ piece of \(W\) is a Seifert manifold. 
A  closed  graph manifold \(W\) is called an \emph{ordinary graph manifold} if:
\begin{enumerate}
\item each Seifert piece \(M_i\) of \(W\)  admits a unique Seifert fibering up to isotopy;
\item the base orbifold of \(M_i\)  is orientable, and \(M_i\) does not admit any singular fiber with slope-\(\frac{q}{2}\) where \(q\) is coprime to \(2\);
\item \(W\) is not homeomorphic to \(M(0,0; \frac{q_1}{p_1}, \frac{q_2}{p_2}, \frac{q_3}{p_3})\) (\(p_1 , p_2, p_3 > 1\)).
\end{enumerate}

It is a classical result that  most Seifert manifolds  admit a unique Seifert fibering up to isotopy \cite[Corollary 3.12]{JWW}. 
From the  definition, ordinary graph manifolds are the majority in  the set of  closed  graph manifolds. In this paper,  we will
only study the indexed links of 
NMS flows on ordinary graph manifolds.
Our results will show that the indexed links of NMS flows on this class of  manifolds 
 are strongly and cleanly related to the topology of the underlying graph manifolds.

\subsubsection{The indexed link   related to the JSJ decomposition}
Let $W$ be an  ordinary graph manifold and  \(W=M_1 \cup \cdots \cup M_s\) be a  JSJ decomposition with the JSJ tori set \(\mathcal{T}\). 
We denote by $W|T$ the manifold obtained  by cutting $W$ along $T\in \mathcal{T}$.
 An indexed link $l$ in $W$ \emph{is related to the JSJ decomposition $W=M_1\cup \cdots \cup M_s$} if: 
\begin{enumerate}
\item \(l \cap (\cup_{T \in \mathcal{T}} T)=\varnothing\) and \(l\) contains both index-\(0\) knots and  index-\(2\) knots.
If \(T\in \mathcal{T}\) is  separating in \(W\), then there is a knot of \(l\) with  index $0$ or  $2$
in 
each connected component of \( W| T\).

\item    For each \(i=1, \cdots ,s\), there is a Seifert fibering of $M_i$, such that
$l_i=l \cap M_i$ 
is a union of fibers which includes all of the singular fibers, and every singular fiber knot  is either index-$0$ or index-$2$. 
 
\item  
Let $x_i$ be the number of index-\(1\) knots of $l_i$, \(z_i\) be the number of other knots of $l_i$, \(b_i\) be the number of boundary components of \(M_i\), 
 and \(g_i\) be the genus of the base orbifold of \(M_i\). Then \( z_i+b_i=x_i-2g_i+2 \).

\end{enumerate}

Notice that the condition (1) in the above definition implies that when $T\in \mathcal{T}$ is separating in $W$, the one connected component of $W|T$ contains at least one index-$0$ knot of $l$ and the other connected component of $W|T$ contains at least one index-$2$ knot of $l$.

 \subsection{Main results}
Our first main result (Theorem \ref{t.realization}) explains that the indexed link    related to a JSJ decomposition of an ordinary graph manifold $W$  can be realized as the indexed link of some NMS flow on $W$.

\begin{theorem}\label{t.realization}
Let $W$ be an ordinary graph manifold and $l$ be an indexed link related to a JSJ decomposition of $W$. Then there exists an NMS flow $\phi_t$ on $W$ such that $l$ is the indexed link of $\phi_t$.

\end{theorem}


The second main result (Theorem \ref{t.main0}) shows  that  the indexed link of an NMS flow on an ordinary graph manifold
can be well understood by 
 Operation of changing regular fibers and  Operation A.  
 Roughly speaking, applying 
Operation of changing regular fibers for an indexed link $l$ is to
select an incompressible torus set of $W$ that  does not intersect with $l$ and
splits $W$ into a finite number of  atoroidal blocks. Then
replace the knots of $l$ in each block homeomorphic to $T^2 \times [0,1]$ with two regular fibers of a Seifert fibering of  this block.
 Note that Operation A  is essentially  consistent with Wada’s operations in \cite{Wa}.
For the precise definitions of
Operation of changing regular fibers and 
 Operation A, we refer to Section \ref{s.change regular fibers} and Subsection \ref{ss.Indlink} respectively.

\begin{theorem}\label{t.main0}
Let \(l\) be an indexed link in an ordinary graph manifold $W$. Then \(l\) is the indexed link 
of an NMS flow on $W$ if and only if
there is an indexed link $l'$ related to 
a JSJ decomposition of $W$ such that $l$
can be obtained by  $l'$ 
 using at most one step of Operation of changing regular fibers and then applying finitely many steps of operations in Operation A. 
\end{theorem}

Notice that applying at most one step of  Operation of changing regular fibers 
 actually means changing a finite number of regular fiber pairs.

\subsection{Further remarks}
It is natural to expect a complete classification   of NMS flows (up to topological equivalence) on a (ordinary) graph manifold.
But generally  heteroclinic trajectories connecting saddle orbits will
lead the question quite wild. As a first step to study this question, a reasonable object is to classify NMS flows without any heteroclinic trajectory.
In \cite{Yu},
the second author of this paper devised a path for the $3$-manifold $S^3$ to discuss this issue and classified such NMS flows completely in orbits with a small number of periods (no more than $4$).  But even on $S^3$, it still seems  difficult to provide a complete classification for such NMS flows.
Nevertheless, it remains an interesting further topic to consider similar problems on (ordinary) graph manifolds.

Readers may wonder how the main results of this paper can be generalized to all graph manifolds. In fact, there is no inherent difficulty, but similar results and proofs would be very complicated and subtle. For specificity and clarity, in this paper we focus on  ordinary graph manifolds.

\subsection{Outline of the article}
 This paper is organized as follows. In Section \ref{s.pre}, we introduce some definitions and elementary properties.
In Section \ref{s.indexed links}, we  discuss the indexed links of NMS flows on ordinary graph manifolds. In Section \ref{ss.Lya}, we  discuss the Lyapunov graphs  of NMS flows. 
In Section \ref{s.JSJ decompositions}, we discuss the indexed link related to a  JSJ decomposition of an ordinary graph manifold, and prove Theorem \ref{t.realization}.
In Section \ref{s.change regular fibers}, we prove  Theorem \ref{t.main0}.

\section{preliminary}\label{s.pre}

\subsection{FRH decompositions}\label{ss.FRHD}

Round handle decomposition was firstly introduced by Asimov in \cite{As1}, which is   closely related to NMS flows. In \cite{Mo}, Morgan expanded Asimov's definition slightly to allow for nonorentability in the stable and unstable manifolds. Now, we introduce the definition of the $3$-dimensional orientable round $k$-handle (abbreviated as  \emph{\(k\)-RH}, $k=0,1,2$), as defined by Morgan.

Let \(X\) be the vector field on \( I \times D^k \times D^{2-k}\)
 given by 

\begin{equation}
\centering
X( \theta , x, y) = \frac{\partial}{\partial \theta} - \sum_{i=1}^k x_i  \frac{\partial}{\partial x_i}+\sum_{j=1}^{2-k} y_j \frac{\partial}{\partial y_j}
\end{equation}
where $I=[0,1]$,  \(x_i\)  and \(y_j\) denote the standard coordinate functions on \(\mathbb{R}^k\) and 
\(\mathbb{R}^{2-k}\) respectively, and \(\frac{\partial}{\partial \theta}\) is a vertor field along the direction of \(I\).

When $k=0,2$, we call the solid torus $I \times D^k \times D^{2-k}/(0, x, y) \sim (1, x, y)$ a \emph{ $k$-RH}.
When $k=1$, we call the solid torus $I \times D^k \times D^{2-k}/(0, x, y) \sim (1, x, y)$ an \emph{untwisted  $1$-RH}, and call the solid torus $I \times D^k \times D^{2-k}/(0, x, y) \sim (1, -x, -y)$ a \emph{twisted  $1$-RH}.
Let $R$ be an RH.
We call the NMS flow on  $R$ induced by $X$ the \emph{natural NMS flow} on $R$. We denote by $\partial_{-}R$ the closure of the incoming boundary set of the natural NMS flow on $R$, and denote by $\partial_{+}R$ the  closure of \(\partial R\smallsetminus \partial_{-}R\).

\begin{figure}[htbp]
\centering
\subfigure[untwisted]{\includegraphics[width=0.25\textwidth]{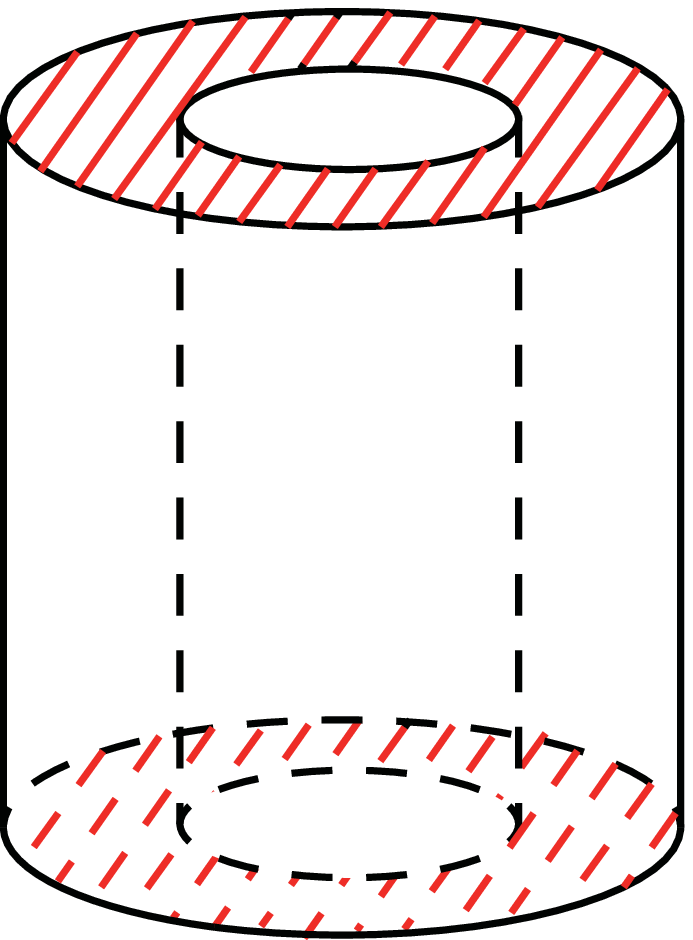}}
\hspace{.90in}
\subfigure[twisted]{\includegraphics[width=0.32\textwidth]{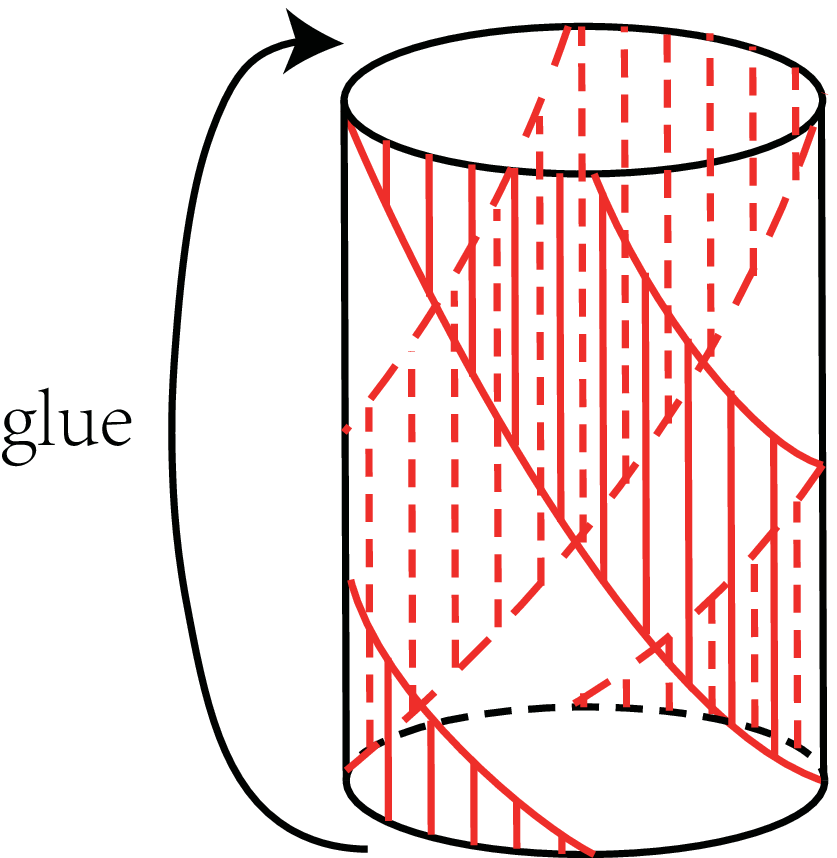}}

\caption{The  shaded part is $\partial_{-}R$.}
\label{Fg6}
\end{figure}

Suppose that  $R$ is a  $1$-RH. If $R$ is untwisted, then $\partial_{-}R$  is a  union of two disjoint annuli (see Figure \ref{Fg6} (a)). If  $R$ is twisted, then $\partial_{-}R$ is an annulus (see Figure \ref{Fg6} (b)). At this point, we call the cores of the connected components of  \(\partial_{-}R\) the  \emph{attaching circles}.

For a compact 3-manifold \(N\), if there exists a non-singular flow on \(N\) pointing inward on  \(\partial _{-}N\) (\(\partial _{-}N\) is a specified union of components of  \(\partial N\)) and outward on  \(\partial _{+}N = \partial N \smallsetminus \partial _{-}N\), then we use the notation \((N,\partial _{-}N)\) to denote it. It should be noted that we do not rule out the cases \( \partial N= \varnothing \), \( \partial_{-} N= \varnothing \) or \( \partial_{+} N= \varnothing \). 
\begin{defi}
Let $M$ and $N$ be two orientable $3$-manifolds.
  \(M\) is obtained from \((N, \partial_{-}N)\) by \emph{attaching a \(k\)-RH} \(R\) if  there is an embedding \(\varphi: \partial_{-}R \rightarrow \partial_{+} N\) such that \(M\cong N\cup_{\varphi }R\) (\(k=0,1,2\)).
\end{defi}

\begin{rema}
We specify that   \(\partial _{+}M= \overline{(\partial _{+} N \smallsetminus \varphi(\partial_{-} R)) \cup \partial_{+}R}\), the  closure of $(\partial _{+} N \smallsetminus \varphi(\partial_{-} R)) \cup \partial_{+}R$.
\end{rema}

\begin{defi}
An orientable \(3\)-manifold \(M\) admits a \emph{round handle decomposition} (abbreviated as \emph{RH decomposition}) \(M=(\partial _{-}M \times I)\cup_{i=1}^{n} R_{i}\) if  each \(R_{i}\) is an RH attached to \((\partial _{-}M \times I)\cup_{j < i} R_{j}\), where \(\partial _{-}M \) is a specified union of components of  \(\partial M\),   and we specify that \(\partial _{+} (\partial_{-} M \times I) = \partial_{-} M \times \{1\}\).
\end{defi}

\par Suppose that \(M\) is obtained from \((N, \partial_{-}N)\) by attaching a \(1\)-RH \(h\), and that both \(M\) and \(N\) are  orientable. 
 Let \(A\) be a union of  components of \(\partial_{+} N\) such that \(\partial_{-}h\) is attached to \(A\) and meets each component of \(A\).  We fatten up \(A\) to get \(A \times I\) such that \(\partial_{-} h\) is attached to \(A \times \{1\}\). 
The manifold \(C(h)\) obtained from \(A \times I\) by attaching   \(h\)  is called  a \emph{fat round \(1\)-handle} (abbreviated as  \emph{\(1\)-FRH}), and we specify that \(\partial _{-} (C(h))= A \times \{0\}\).
 Making \(\partial _{-} (C(h))\) naturally attach to  \( \partial_{+} N\), we still get \(M\). Namely, \(M\) is obtained  from \(N\) by attaching a \(1\)-FRH \(C(h)\).
Therefore for an RH decomposition of \(M\), we can get a new decomposition of \(M\) by replacing each \(1\)-RH with the corresponding \(1\)-FRH, and call this new decomposition a  \emph{fat round handle decomposition}  (abbreviated as \emph{FRH decomposition}). In addition, we also refer to \(0\)-RHs as \emph{\(0\)-FRHs}, and refer to \(2\)-RHs as \emph{\(2\)-FRHs}.

Morgan \cite[Lemma 3.1]{Mo} proved the following result that will be very useful in this paper. 

\begin{lemm}
\label{lem1}
Let \(W\) be a graph manifold, and \(W=(\partial _{-}W \times I)\cup_{i=1}^{n} R_{i}\) be an RH decomposition of \(W\). Let  \(W_{j}= (\partial _{-}W \times I)\cup_{i=1}^{j} R_{i}\) for \(j=1,\cdots, n\), then every \(\partial W_{j}\) is a  union of pairwise disjoint tori.
\end{lemm}

In  \cite{Wa},
Wada proved that the $1$-FRH realized in FRH decompositions of $S^3$ is one of the types (a)-(e) in  Lemma \ref{lem2}. Let \(W\) be a graph manifold, $h$ be a $1$-RH in an RH decomposition of $W$, and 
\(C(h)\) be the \(1\)-FRH associated to \(h\).
 In Proposition 3.5 of \cite{Mo}, Morgan proved that   \(C(h)\)  is one of types (a), (b), (c) and (g) in  Lemma \ref{lem2} when there is an inessential attaching circle in \(C(h)\).
In fact, by using both the tools and the results due to Morgan \cite{Mo} and Wada   \cite{Wa},  it is not difficult to get all of the $1$-FRHs realized in FRH decompositions of $W$. The following lemma describes the complete classification in question.  Note that the proof is only a slight generalization of the tools in Morgan \cite{Mo} and Wada   \cite{Wa},  so we only give the sketch of a proof here.

\begin{lemm}
\label{lem2} 
 Let  \(W\) be a graph manifold, $h$ be a $1$-RH in an RH decomposition of $W$, and 
\(C(h)\) be the \(1\)-FRH associated to \(h\). Let \(r\) be a core of \(h\). Then
 \((C(h),r)\) is one of the following types:
\begin{enumerate}[\indent (a)]

\item \(C(h)\cong (T_{1}\times I) \sharp (T_{2}\times I)\), that is, the connected sum of  \(T_{1}\times I\) and \(T_{2}\times I\), where \(T_{1}\) and \(T_{2}\) are tori. Moreover,  \(\partial _{-}(C(h))=(T_{1}\times  \left \{ 0 \right \}) \cup (T_{2}\times  \left \{ 0 \right \})\), and \(r\) bounds a disk in \(C(h)\).

\item \(C(h)\cong( T^{2}\times I) \sharp (D^{2} \times S^{1})\), \(\partial _{-}(C(h))=T^{2}\times \left \{   0\right \}\) or \(\partial _{-}(C(h))=(T^{2}\times \left \{ 0  \right \}) \cup (\partial D^{2}\times S^{1})\), and \(r\) bounds a disk in \(C(h)\).
\item \(C(h)\cong V_{1} \sharp V_{2}\) where \(V_{1}\) and \(V_{2}\) are two solid tori, \(\partial _{-}(C(h))= \partial V_{1}\), and \(r\) bounds a disk in \(C(h)\).
\item \(C(h)\cong F \times S^{1}\) where \(F\) is a disk with two holes. Moreover, \(\partial _{-}(C(h))\) is a (connected) component or a union of two components of \(\partial C(h)\), and \(r= \left \{   *\right \} \times S^{1}\) for some point \(*\) in \({\rm Int}F\), i.e. the interior of $F$.
\item \(C(h)\cong D^{2}\times S^{1} \smallsetminus {\rm Int}N\) where \(N\) is a tubular neighborhood of the (2,1)-cable of \(\left \{ 0 \right \}\times S^{1}\) in \(D^{2}\times S^{1}\), \( \partial _{-}(C(h))= \partial N\), and \(r=\left \{ 0 \right \} \times S^{1}\).

\item \(C(h)\cong P \widetilde{\times }S^{1}\) where \(P\) is a  Möbius
strip  with one hole,  \(\partial _{-}(C(h))\) is a component of \( \partial C(h)\), and \(r= \left \{   *\right \} \widetilde{\times} S^{1}\) for some point \(*\) in \({\rm Int} P\).

\item \(C(h)\cong \mathbb{R}\mathbb{P}^{3} \sharp (T^{2} \times I)\),  \(\partial _{-}(C(h))= T^{2} \times \left\{0 \right \}\), and  a \((2,1)\)-cable of \(r\)  bounds a disk in \(C(h)\). Moreover, \(C(h)\cong \mathbb{R}\mathbb{P}^{3} \sharp (T^{2} \times I)\) exists only if \(W \cong \mathbb{R}\mathbb{P}^{3} \).
\end{enumerate}
\end{lemm}
\begin{proo}[Sketch of a proof]
Suppose \(C(h) \cong A \times I \cup_{\varphi}h\)  where \(\varphi: \partial_{-} h  \rightarrow A \times \{1\}\).  By Lemma  \ref{lem1}, \(A\) is a union of  pairwise disjoint tori.

Case 1. \(h\) is untwisted, then \(\partial_{-}h\) consists of two annuli, which implies that \(A\) is either a torus or a  union of two disjoint tori. Let \(c_1\) and \(c_2\) be two attaching circles that are contained in 
different components of \(\partial_{-}h\).
 We endow  the orientations on \(c_1\) and \(c_2\) such that $c_1$ and $c_2$ are isotopic in $h$.

Case 1.1. \(A\) is the   union of two disjoint tori.
\begin{itemize}
\item If both \(\varphi(c_1)\) and \(\varphi(c_2)\) are essential in \(A \times \{1\}\), then we get type (d).
\item If one of \(\varphi(c_1)\) and \(\varphi(c_2)\) is essential in \(A \times \{1\}\) and the other is inessential in \(A \times \{1\}\), then we get type (b).
\item If both \(\varphi(c_1)\) and \(\varphi(c_2)\) are inessential in \(A \times \{1\}\), then either we get a manifold with a \(S^2\) boundary component (see Figure \ref{FgFRH1} (a)) or we get type (a) (see Figure \ref{FgFRH1} (b)). By Lemma \ref{lem1}, the first case is impossible.

\end{itemize}
\begin{figure}[htbp]
\centering
\subfigure[ ]{\includegraphics[width=0.46\textwidth]{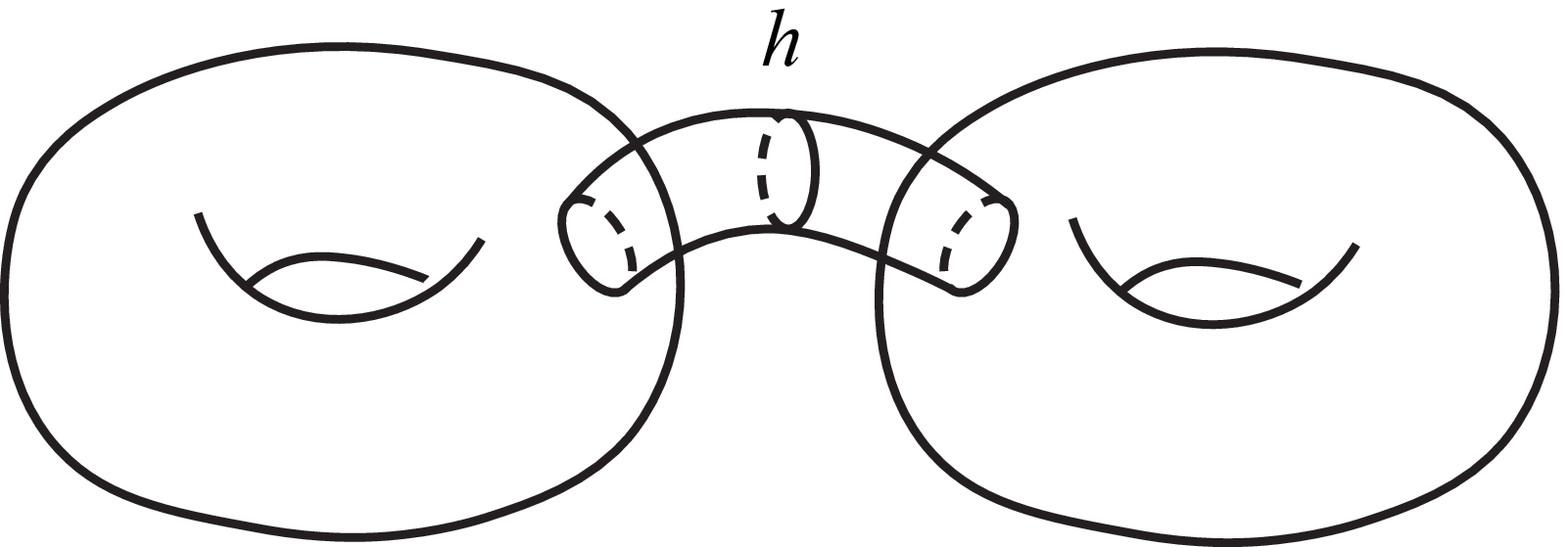}}
\hspace{.30in}
\subfigure[ ]{\includegraphics[width=0.46\textwidth]{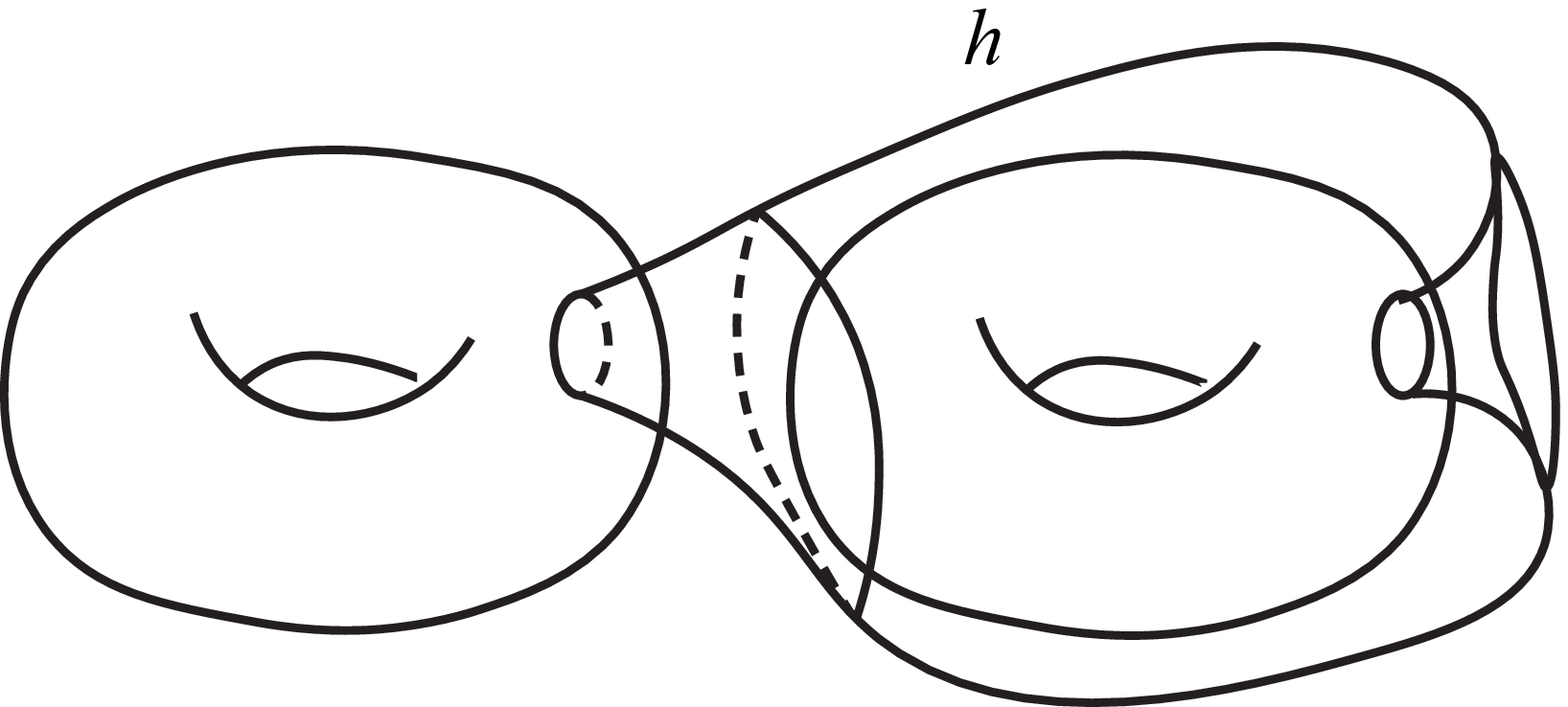}}
\subfigure[ ]{\includegraphics[width=0.52\textwidth]{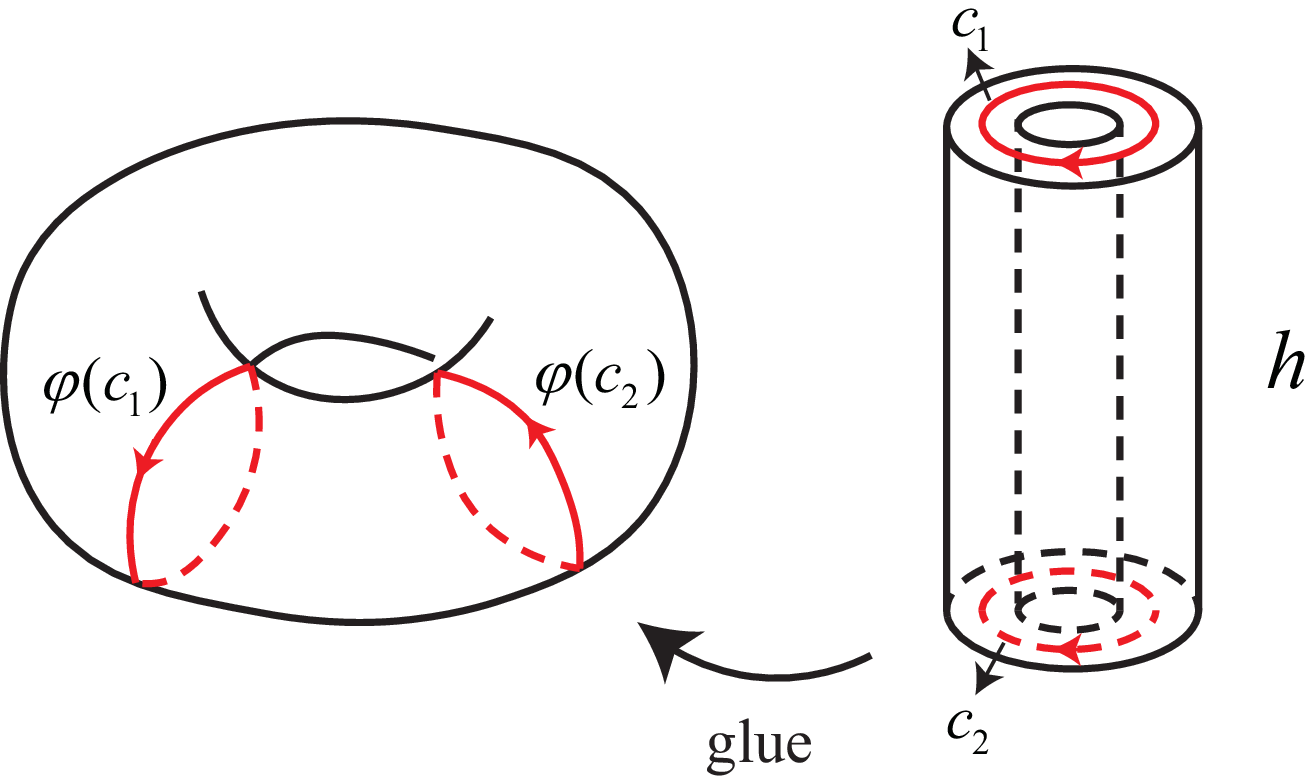}}
\caption{ }
\label{FgFRH1}
\end{figure}

Case 1.2. \(A\) is a torus.

\begin{itemize}
\item Suppose both \(\varphi(c_1)\) and \(\varphi(c_2)\) are essential in \(A \times \{1\}\). If \(\varphi(c_1)\) and \(\varphi(c_2)\) are in opposite directions on \(A \times \{1\}\), then we get type (f) (see Figure \ref{FgFRH1} (c)). Otherwise, we get type (d).
\item If one of \(\varphi(c_1)\) and \(\varphi(c_2)\) is essential in \(A \times \{1\}\) and the other is inessential in \(A \times \{1\}\), then we get type (c).
\item If both \(\varphi(c_1)\) and \(\varphi(c_2\)) are inessential in \(A \times \{1\}\), then there are four cases as shown in Figure \ref{FgFRH_2}. The first manifold  admits a \(S^2\) boundary component, and both the third and fourth manifolds  contain  non-separating \(2\)-spheres. By Lemma \ref{lem1} and the irreducibility of \(W\), these three cases are impossible. The second manifold is of type (b).

\end{itemize}

 Case 2. \(h\) is twisted, then \(\partial_{-}h\) is an annulus, which implies that \(A\) is  a torus. Let \(c_1\) be an attaching circle.
\begin{itemize}
\item If \(\varphi(c_1)\) is  essential in \(A \times \{1\}\), then we get type (e).
\item Otherwise, we get type (g).

\end{itemize}

\begin{figure}[htbp]
\centering
\subfigure[ ]{\includegraphics[width=0.34\textwidth]{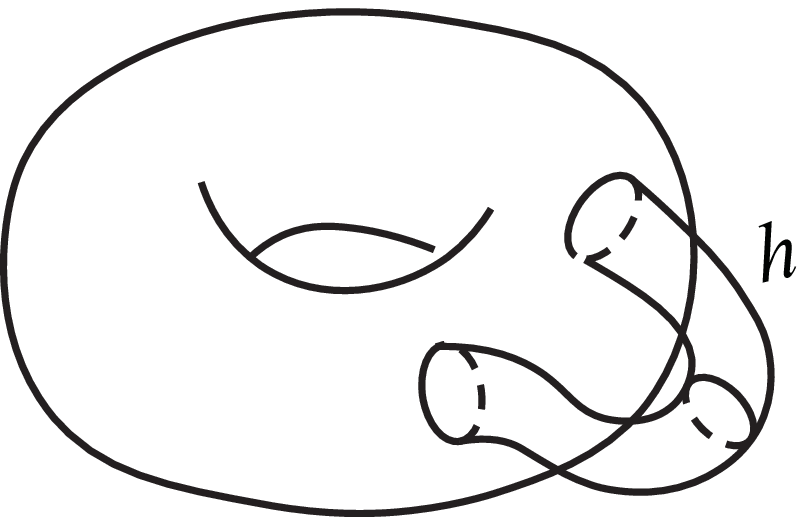}}
\hspace{.60in}
\subfigure[ ]{\includegraphics[width=0.36\textwidth]{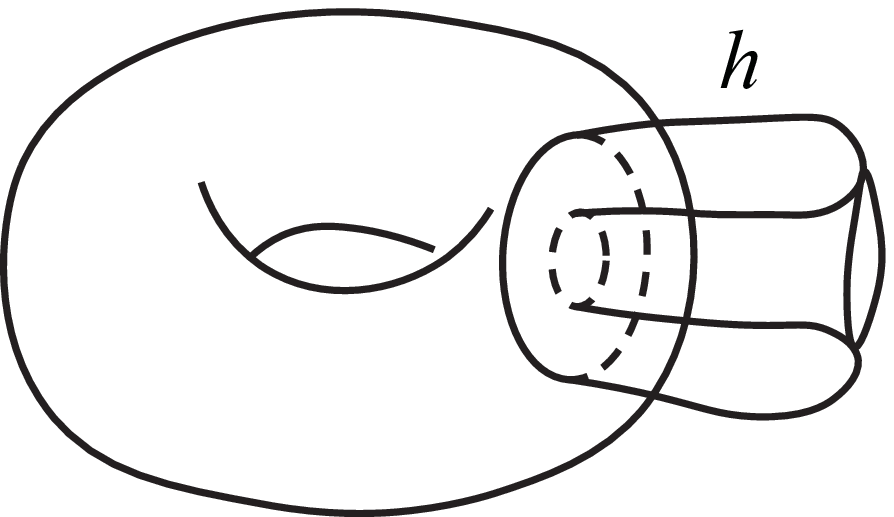}}
\hspace{.60in}
\subfigure[ ]{\includegraphics[width=0.34\textwidth]{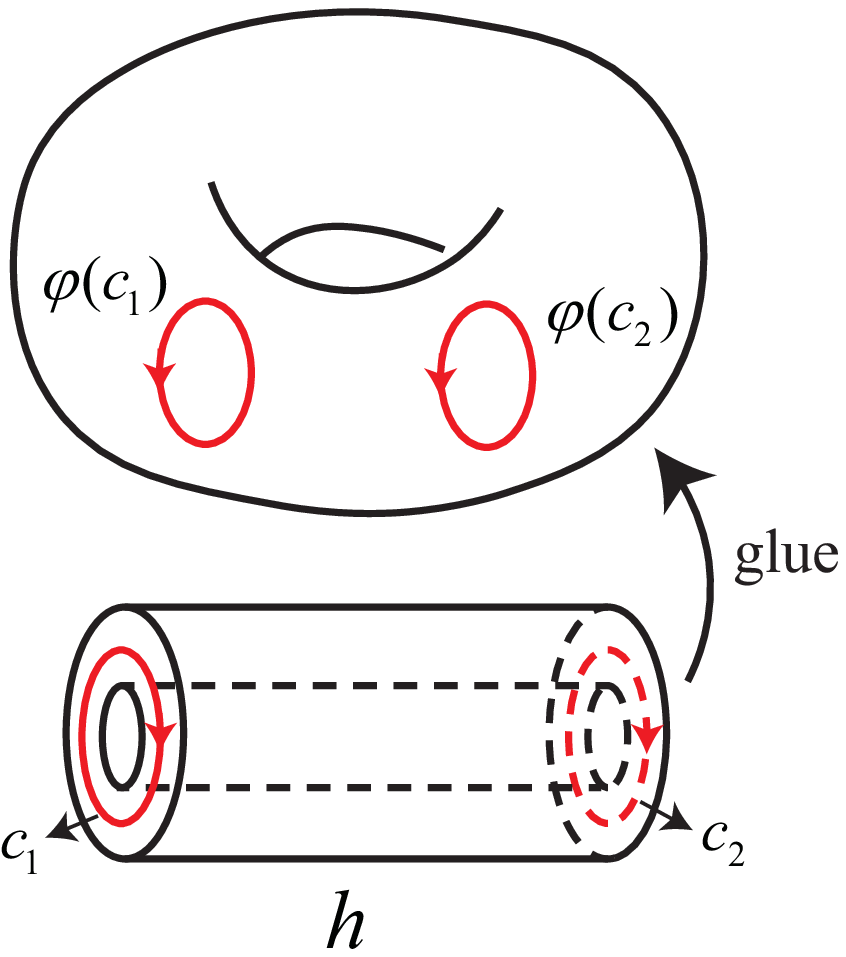}}
\hspace{.60in}
\subfigure[ ]{\includegraphics[width=0.36\textwidth]{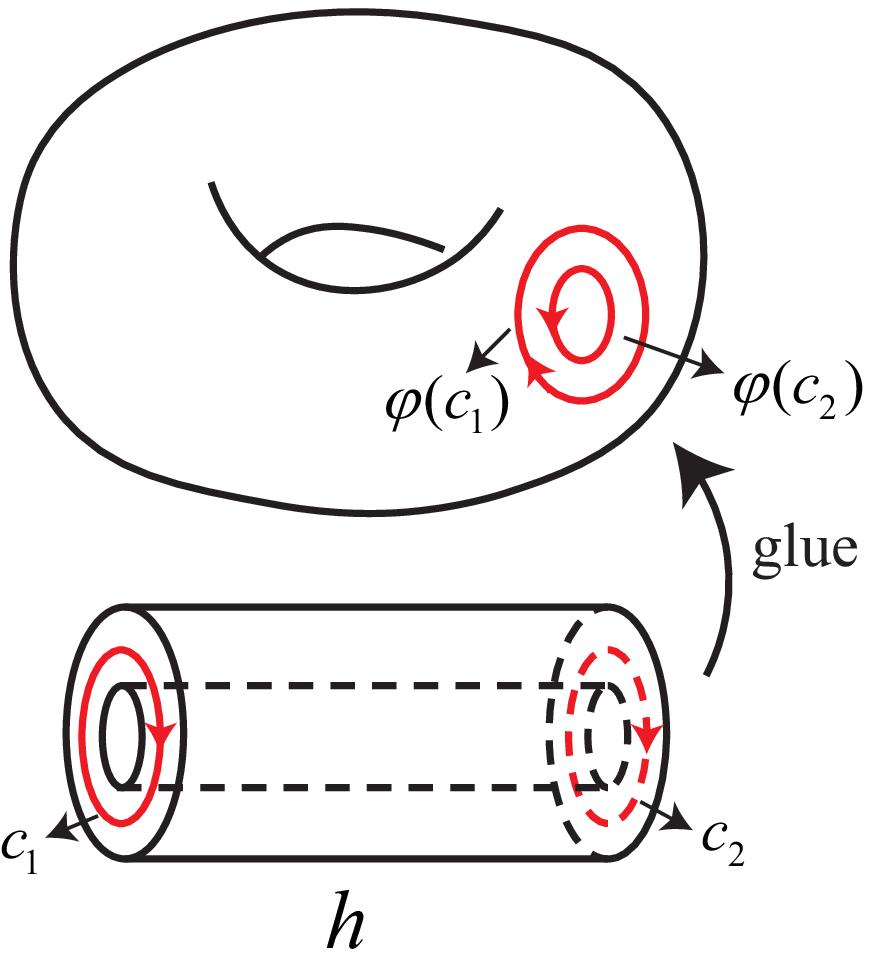}}

\caption{ }
\label{FgFRH_2}
\end{figure}

\end{proo}

\begin{defi}
\label{def3.4}
Let \(C(h)\) be a \(1\)-FRH 
 in a FRH decomposition of a \(3\)-manifold.
We define \(\widetilde{C}(h)\) to be the union of  \(C(h)\) and all of the \(0\)-FRHs and \(2\)-FRHs adjacent to \(C(h)\) in this FRH decomposition.
\end{defi}

 \(\widetilde{C}(h)\) was firstly introduced by Morgan in \cite{Mo} to describe the local situation of the \(1\)-FRH. 
Let \(C(h)\) be a \(1\)-FRH in a FRH decomposition of a graph manifold $W$.
Suppose that \(C(h)\) is 
of type (d), (e) or (f) and   \(\partial \widetilde{C}(h) \neq \varnothing\).
 In Proposition 3.5 of \cite{Mo}, Morgan  got \(\widetilde{C}(h)\) is one of types (4), (5) and (7) in Lemma \ref{lem3}, if there is an inessential attaching circle  in \(\widetilde{C}(h)\). Moreover, he found that  in other cases \(\widetilde{C}(h)\) is a Seifert manifold with a Seifert fibering extended from the attaching circles.  In the following lemma, we describe all cases of \(\widetilde{C}(h)\) in detail.

\begin{lemm} 
\label{lem3}
Let \(C(h)\) be a \(1\)-FRH in a FRH decomposition of a graph manifold $W$.
Suppose that \(C(h)\) is 
of type (d), (e) or (f) and   \(\partial \widetilde{C}(h) \neq \varnothing\).  Then \(\widetilde{C}(h)\)  is one of the following types:
\begin{enumerate}[\indent (1)]
\item \(\widetilde{C}(h) \cong C(h)\);
\item \(\widetilde{C}(h) \cong M (0,2;\frac{q}{p})\);
\item \(\widetilde{C}(h) \cong M (0,1;\frac{q_{1}}{p_{1}},\frac{q_{2}}{p_{2}} )\), where \(p_{1}, p_{2} > 1\);
\item \(\widetilde{C}(h) \cong S^{1} \times D^{2}\);
\item \(\widetilde{C}(h) \cong (S^{1} \times D^{2}) \sharp (S^{1} \times D^{2})\);
\item \(\widetilde{C}(h) \cong M (-1,1;\frac{q}{p})\);
\item \(\widetilde{C}(h) \cong W \sharp (S^{1} \times D^{2})\). 
\end{enumerate}
\par In particular, \(\widetilde{C}(h)\) of type (7) occurs only when \(W\) is homeomorphic to a lens space.
\end{lemm}
\begin{proo} 
If there is no \(0\)-FRH or \(2\)-FRH adjacent to \(C(h)\), then \(\widetilde{C}(h) \cong C(h)\). In the  following of the proof, we consider the case that \(C(h)\) is adjacent to some \(0\)-FRHs or \(2\)-FRHs.
\par Case 1. \(C(h)\) is of type (d).

The three boundary components of \(C(h)\) are symmetric in topology, so we only need to consider two cases:  attaching one solid torus  to \(C(h)\) or attaching two solid tori to \(C(h)\).  Let \(T_i\) be the boundary component of \(C(h)\) for \(i=1,2,3\).
There is a circle bundle \(\pi: F\times S^{1} \cong C(h)\rightarrow F\), where \(F\) is a pair-of-pants, i.e., a disk with two holes. We choose a cross section \(s:  F\rightarrow F\times S^{1}\), endow \(C(h)\) with an orientation, and fix an orientation on \(\partial C(h)\). Let \(d_{i}=s \cap T_i\), and \(l_{i}\) be a fiber in \(T_{i}\).

\par Case 1.1. Attach one solid torus  to \(C(h)\). 

 We may assume that a solid torus \(R_1\) is attached to \(C(h)\) along \(T_1\). let \(\widetilde{m_{1}}\) be a meridian of \(\partial R_{1}\). Suppose that  \(\widetilde{C}(h) \cong C(h)\cup_{\psi}  R_1\), where \(\psi : \partial R_1 \to T_1\) is a diffeomorphism such that \(\psi (\widetilde{m_{1}} )=p_{1}d_{1}+q_{1}l_{1}\) (\(p_1 \geq 0\)). If \(p_{1}=0\), then \(\widetilde{C}(h) \cong (S^{1} \times D^{2}) \sharp (S^{1} \times D^{2})\). Otherwise,  \(\widetilde{C}(h) \cong M (0,2;\frac{q_{1}}{p_{1}})\). In particular,  \(\widetilde{C}(h) \cong T^2 \times I\), if \(p=1\).
\par Case 1.2. Attach two solid tori to \(C(h)\).
\par We may assume that two solid tori \(R_1\) and \(R_2\) are attached to \(C(h)\) along \(T_1\) and \(T_2\) respectively. Let \(\widetilde{m_{j}}\) be a meridian of \(\partial R_{j}\) (\(j=1,2\)).  Suppose that \(\widetilde{C}(h) \cong R_1 \cup_{\psi_{1}} C(h) \cup_{\psi_{2}} R_2\), 
where \(\psi_{j} : \partial R_j \to T_j\) is a diffeomorphism such that \(\psi_{j} (\widetilde{m_{j}} )=p_{j}d_{j}+q_{j}l_{j}\) (\(p_j \geq 0\)). 

If \(p_{1} = p_{2}=0\), then there is a non-separating \(S^{2}\) in \(W\), which contradicts to the irreducibility of \(W\). 
 If one of \(p_{1}\) and \( p_{2}\) is \(0\) and
the other is not, then
 we may assume that \(p_{1}=0\) and \(p_{2} \neq 0\). Then \(C(h) \cup_{\psi_1} R_1\cong V_{1} \sharp  V_{2}\), where \(V_1\) and \(V_2\) are two solid tori,  and \(\partial V_2=T_2 \). Thus  \(\widetilde{C}(h) \cong V_{1} \sharp \left( V_{2}\cup _{\psi_{2}  }R_{2}\right) \). If \(W\) isn't homeomorphic to any lens space, then \(V_{2}\cup _{\psi_{2} }R_{2} \cong S^{3}\), which implies that \(\widetilde{C}(h) \cong S^{1} \times D^{2}\). Otherwise, \(V_{2}\cup _{\psi_{2} }R_{2} \cong W\) or \(S^{3}\), then \(\widetilde{C}(h) \cong W \sharp (S^{1} \times D^{2})\) or \(S^{1}\times D^{2}\). 

If  \(p_{1}=1\) or \( p_{2}=1\)  and we assume that \(p_{1}=1\) without confusion, then \(C(h) \cup_{\psi_1} R_1 = T^2 \times I\), which implies that \(\widetilde{C}(h) \cong S^1 \times D^2\).
 If \(p_1, p_2 \geq 2\), then \(\widetilde{C}(h) \cong M (0,1;\frac{q_{1}}{p_{1}},\frac{q_{2}}{p_{2}})\).

\par Case 2.  \(C(h)\) is of type (e), then \(C(h) \cong M(0,2; \frac{1}{2})\). 
\par  Removing a small open tubular neighborhood \(N\) of the singular fiber of \(C(h)\), we get a manifold \(F \times S^1\) with a circle bundle \(\pi: F \times S^1  \rightarrow F\), where \(F\) is a pair-of-pants.
Similar to Case 1, it is easy to prove that \(\widetilde{C}(h)\) is one of types (3), (4) and (7). In particular, if  \(\widetilde{C}(h)\) is of type (7), then \(W \cong \mathbb{R} \mathbb{P}^3\). If \(\widetilde{C}(h)\) is a Seifert manifold, then \(\widetilde{C}(h)\) admits 
 a Seifert fibering  such that 
 \(\widetilde{C}(h)\)  contains a singular fiber with slope-\(\frac{q}{2}\), where \(q\) is coprime to \(2\).
\par Case 3. \(C(h)\) is of type (f), then \(C(h) \cong M(-1,2;)\).

Similarly, it is easy to prove that \(\widetilde{C}(h)\) is of type (6). 

\end{proo}

\begin{rema}\label{paichu}
Suppose that \(\widetilde{C}(h)\) is a Seifert manifold with boundary, then: 
\begin{enumerate}
\item
if \(C(h)\) is of type (e), then   \(\widetilde{C}(h)\) admits 
 a Seifert fibering  such that 
 \(\widetilde{C}(h)\)  contains a singular fiber with slope-\(\frac{q}{2}\), where \(q\) is coprime to \(2\);
\item
if \(C(h)\) is of type (f), then  \(\widetilde{C}(h)\) admits a Seifert fibering such that the base orbifold of \(\widetilde{C}(h)\) is non-orientable.
\end{enumerate}

\end{rema}

 In the following, when we talk about a flow on a manifold with boundary, we mean that the flow  is transverse to the boundary of this manifold.
Let \(\phi_t\) be an NMS flow on an orientable \(3\)-manifold \(M\), and $c$ be a closed orbit   of $\phi_t$.
A \emph{filtrating neighborhood} $N$ of    $c$   is a connected neighborhood of $c$ with boundary such that 
  $\phi_t$ is transverse to   $\partial N$ and
  the maximal invariant set of $\phi_t|_{N}$   is  $c$. 

Let \(M =(\partial _{-}M \times I) \cup_{i=1}^{n} C_i\) be a FRH decomposition of  $M$   associated to an RH decomposition $M =(\partial _{-}M \times I) \cup_{i=1}^{n} R_i $, where $C_i$ is a FRH associated to the RH $R_i$. According to the definition of FRH decompositions,
if  $C_i$ is a $1$-FRH, then $C_i $ is obtained from $R_i$ by attaching some thickened surfaces. Otherwise $C_i$ is the same as $R_i$.
%
%
%

 Suppose that $\phi_t$ has $n$ closed orbits exactly,
  is   transverse 
inwardly to $\partial _{-}M$, and is   transverse 
outwardly to $\partial M \smallsetminus\partial _{-}M$.
 We say that the   decomposition \(M = (\partial_{-}M \times I) \cup_{i=1}^{n} C_i\)  is \emph{a FRH decomposition of \(\phi_t\)}  if 
\begin{itemize}
\item $C_i$ is a filtrating neighborhood of a closed orbit $c_i$ of $\phi_t$;
\item
the components of $\partial C_i$ that are transversed inwardly by $\phi_t|_{C_i}$ is the gluing region $\partial_{-}C_i$  exactly;

\item $c_i$ is a core of the RH $R_i$ associated to the FRH $C_i$.

\end{itemize}

Suppose that  \(M = (\partial_{-}M \times I) \cup_{i=1}^{n} C_i\)  is a FRH decomposition of \(\phi_t\). It is easy to observe that $C_i$ is a $k$-FRH if and only if the index of the closed orbit $c_i$  is $2-k$ for $k=0,1,2$.
FRH decompositions and NMS flows are closely related:

\begin{theorem}
\label{t.thm1} 
Let $M$ be an orientable $3$-manifold.
If \(M\) admits a  FRH decomposition, then there exists an NMS flow \(\phi_t\) on \(M\) such that this decomposition is a  FRH decomposition of \(\phi_t\).
Conversely, if \(M\) admits an NMS flow \(\phi_t\), then there exists a  FRH decomposition of \(\phi_t\).

\end{theorem}
Theorem \ref{t.thm1}  is obtained entirely by the proof of Theorem N of Asimov \cite{As1} and Proposition in page 43 of Morgan \cite{Mo}, so we have omitted the proof here.

%
%
%
%
%
%
%
%

Let $M$ be an orientable $3$-manifold admitting an NMS flow $\phi_t$.
By Theorem \ref{t.thm1},  the FRH decomposition of $\phi_t$  always exists.
We can construct a FRH decomposition of $\phi_t$ using a  Lyapunov function $f$ associated to $\phi_t$, as follows.
We  choose finitely many  regular level sets of $f$ to decompose $(M, \phi_t)$ into some filtrating neighborhoods. By Theorem 4.4 in Yu \cite{Y}, it is easy to prove that every filtrating neighborhood 
  is a FRH.  Thus we construct
a FRH decomposition
   of \(\phi_t\).
Different Lyapunov functions associated to $\phi_t$ maybe induce different FRH decompositions of $\phi_t$. Thus, 
 the FRH decompositions of $\phi_t$  may  not be unique.


\subsection{Indexed links}\label{ss.Indlink}
\par From now on, we define some operations to discuss the indexed link of an NMS flow on a  graph manifold \(W\). Let \(k\) be an indexed knot in \(W\). We denote by \(N(k, W)\) the regular neighborhood of \(k\) in \(W\), and denote by \({\rm Ind}(k)\) the index of \(k\).

\par Let \(l'\) be a link in \(W\),
and \(l''\) be a link in \(S^3\). Let \(N'\) (resp. \(N''\)) be a small open tubular neighborhood of \(l'\) (resp. \(l''\))  in \(W\) (resp. \(S^3\)). Delete  small open \(3\)-balls \(B'\) and \(B''\) from  \(W\smallsetminus N'\) and \(S^3\smallsetminus N''\) respectively, and we obtain a new link in \(W\) by gluing \(W\smallsetminus B'\) and \(S^3\smallsetminus B''\) along their \(2\)-sphere boundaries. We denote by \(l' \cdot l''\)
this new link in \(W\).

\begin{figure}[htbp]
\centering
\includegraphics[scale=0.3]{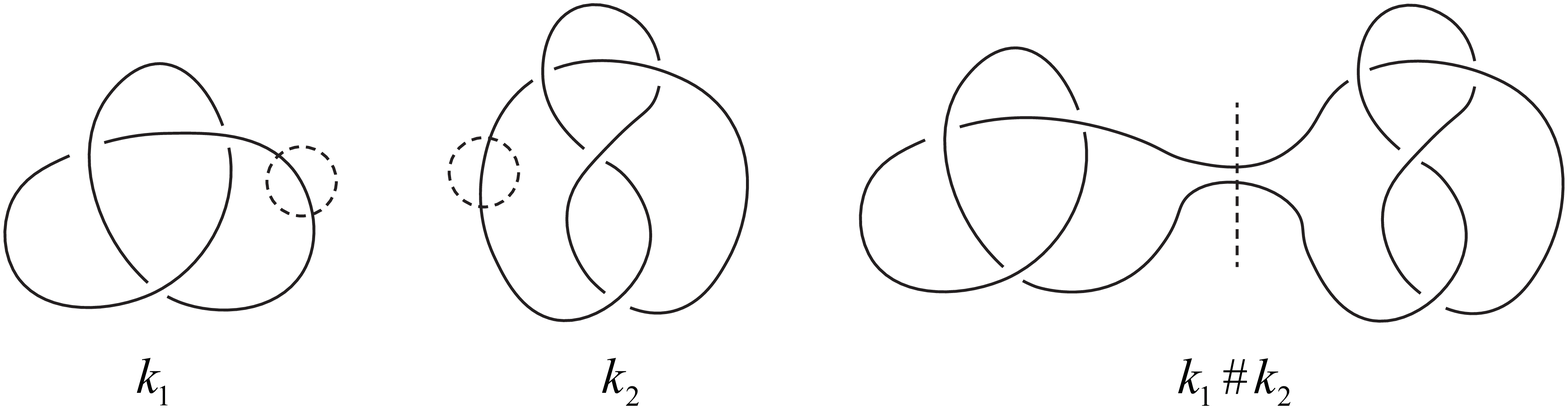}
\caption{ }
\label{Fgconsum}
\end{figure}

Let \(k_1\) be a knot in \(W\),
and \(k_2\) be a knot in \(S^3\). We choose an open \(3\)-ball \(B_1\) (resp. \(B_2\)) in \(W\) (resp. \(S^3\)) such that \(B_1 \cap k_1\) (resp.  \(B_2 \cap k_2\)) can be embedded in a \(2\)-sphere.
By gluing \(W\smallsetminus B_1\) and \(S^3\smallsetminus B_2\) along their  \(2\)-sphere boundaries,
we can  get a new knot  in \(W\), called the \emph{connected sum} \(k_1 \sharp k_2\)  
(see Figure \ref{Fgconsum}). 

Both \(l'\) and \(l''\) are defined as before, and let 
 \(k'\) (resp. \(k''\)) be a knot of \(l'\) (resp. \(l''\)).  The \emph{connected sum \(l' \sharp l''\) associated to \(k' \sharp k''\)} is a new link in \(W\) obtained from \(l'\) and \(l''\) by doing connected sum \(k' \sharp k''\). Obviously, \(l' \cdot l''\) and  \(l'_1 \sharp l'_2\) associated to \(k' \sharp k''\) can be realized by doing connected sum \(W \sharp S^3\).

The terminologies of knot theory refer to Rolfsen \cite{Rd}. 
 For a given  indexed link \(l_1\) in a graph manifold \(W\), we define seven operations as follows.
\par
\begin{figure}[htbp]
\centering
\includegraphics[scale=0.42]{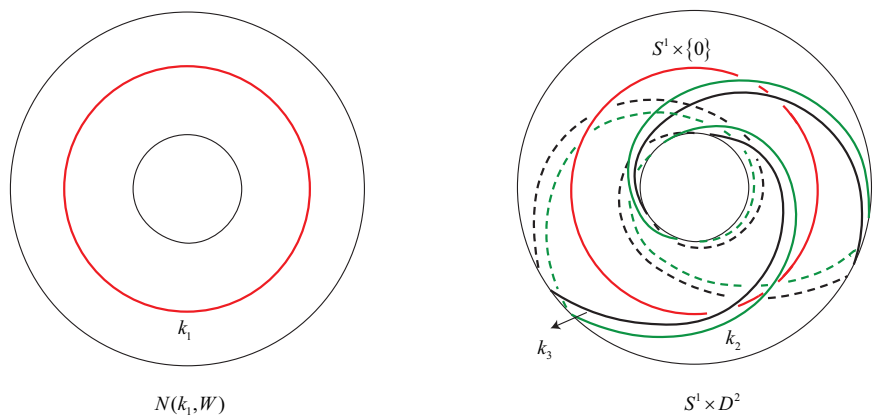}
\caption{\(k_{2}\) and \(k_{3}\) are two parallel \(\left ( 3,2 \right )\)-cables of \(S^{1}\times \left \{ 0 \right \}\).}
\label{FgOA1}
\end{figure}

\begin{figure}[htbp]
\centering
\includegraphics[scale=0.42]{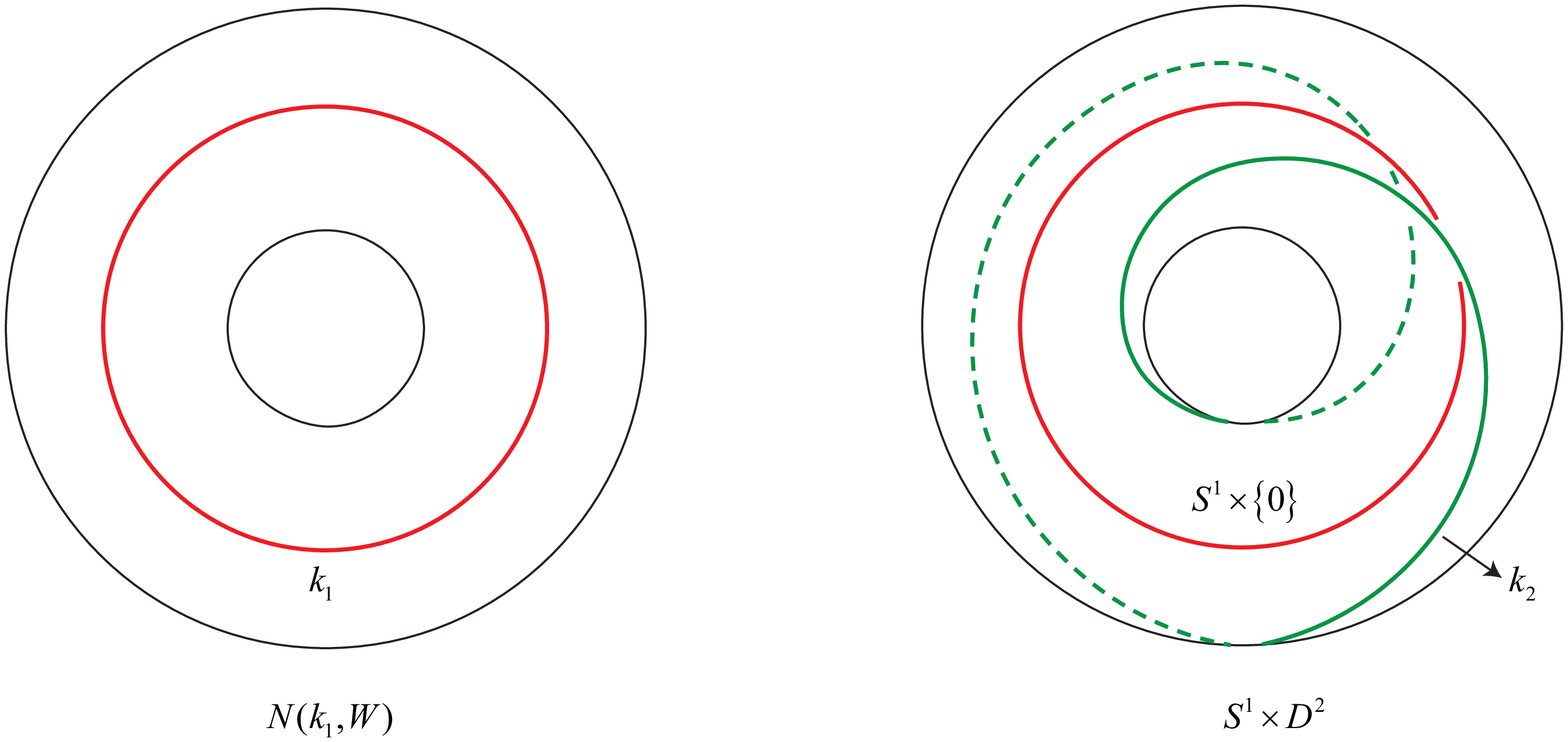}
\caption{\(k_{2}\) is a \(\left ( 2,1 \right )\)-cable of \(S^{1}\times \left \{ 0 \right \}\).}
\label{FgOA2}
\end{figure}

 \textbf{Operation A:} Choose an indexed link \(l_2\) of an NMS flow on \(S^{3}\) and  an unknot  \(u\) in \(S^3\) with index 1.
\begin{enumerate}[\indent I.]
\item To make \(l_{1}\cdot l_{2}\cdot u\).
\item To make \(l_{1}\cdot (l_{2}\smallsetminus k_{2})\cdot u\), where \({k}_2\) is a knot of \({l}_2\) with index \(0\) or \(2\).

\item To make \((l_{1}\smallsetminus k_{1})\cdot l_{2}\cdot u\), where \({k}_1\) is a knot of \({l}_1\) with index \(0\) or \(2\).

\item To make \((l_{1}\smallsetminus k_{1})\cdot (l_{2}\smallsetminus k_{2})\cdot u\), where  \({k}_i\) is a knot of \({l}_i\) (\(i=1,2\)), and \({\rm Ind}({k}_1)= 2- {\rm Ind}({k}_2)=0\) or \(2\).
\item To make \((l_{1} \sharp l_{2}) \cup m\). The connected sum \(l_{1} \sharp l_{2}\) is associated to
\(k_1 \sharp k_2\), where \(k_i\) is a knot of \(l_i\)  with index \(0\) or \(2\) (\(i=1,2\)). \({\rm Ind}(k_{1} \sharp k_{2}) = {\rm Ind}(k_{1})\) or \({\rm Ind}(k_{2})\), and \(m\) is a meridian of \(k_{1} \sharp k_{2}\) with index \(1\).
\item Choose a knot \(k_{1}\) of \(l_{1}\) with index \(0\) or \(2\), and replace \(N(k_{1},W)\) by \(S^{1} \times D^{2}\) with three indexed circles in it: \(S^{1}\times \left \{ 0 \right \}\), \(k_{2}\) and \(k_{3}\). Here, \(k_{2}\) and \(k_{3}\) are two parallel \(\left ( p,q \right )\)-cables of \(S^{1}\times \left \{ 0 \right \}\) (see Figure \ref{FgOA1}). The indices of \(S^{1}\times \left \{ 0 \right \}\) and \(k_{2}\) are either \(0\) or \(2\), and one of them is equal to \({\rm Ind}(k_{1})\). \({\rm Ind}(k_{3}) =1\).
\item Choose a knot \(k_{1}\) of \(l_{1}\) with index \(0\) or \(2\), and replace \(N(k_{1},W)\) by \(S^{1} \times D^{2}\) with two indexed knots in it: \(S^{1}\times \left \{ 0 \right \}\) and a \(\left ( 2,q \right )\)-cable knot \(k_{2}\) of \(S^{1}\times \left \{ 0 \right \}\) (see Figure \ref{FgOA2}). \({\rm Ind}(S^{1}\times \left \{ 0 \right \}) =1\), and \({\rm Ind}(k_{2})={\rm Ind}(k_{1})\). 
\end{enumerate}

In  \cite{Wa}, Wada 
classified the indexed links of NMS flows on \(S^3\) by using operations similar to Operation A.
In fact, Operation A is consistent with Wada's operations in \cite{Wa}, except that Wada's operations  require that  both \(l_{1}\) and \(l_{2}\) are two given indexed links  in \(S^{3}\).

\subsection{NMS flows related to incompressible torus decompositions}\label{s.ITD}

An NMS flow  \(\phi_t\) on a graph manifold \(W\) is called an \emph{NMS flow related to  incompressible torus decompositions}  if:
\begin{enumerate}
\item the indexed link \(l\) of \(\phi_t\) contains index-\(1\) knots, and \(l\)
is \emph{knotted} in \(W\), i.e., no connected component
of $l$ bounds an embedded disk in $W$;
\item there is a FRH decomposition  \(W=  (\partial_{-}W \times I) \cup_{i} C_i\) of \(\phi_t\) such that  every boundary surface of each 1-FRH  is either  incompressible in \(W\) or  a boundary of some 0-FRH or 2-FRH in this decomposition.
\end{enumerate}

Note that the \(0\)-FRH is the same as the \(0\)-RH, and the \(2\)-FRH  is the same as the \(2\)-RH.
At this point, we also say that $\phi_t$ is related to an incompressible torus decomposition \(W=  (\partial_{-}W \times I) \cup_{i} C_i\), and we can write this decomposition  in the form of $W=  (\partial_{-}W \times I) \cup_{j} \widetilde{C}(h_{j})$ (where $\widetilde{C}(h_{j})$ is defined in Derfintion \ref{def3.4}).


\begin{lemm}
\label{lemITD}

Suppose that \(\phi_t\) is an NMS flow on a  graph manifold \(W\) with saddle closed orbits, and  \(\phi_t\) admits at least two saddle closed orbits when \(\partial W = \varnothing\). Then 
\(\phi_t\) is related to an  incompressible torus decomposition
 \(W=  (\partial_{-}W \times I) \cup_{j=1}^{n}\widetilde{C}(h_{j})\) if and only if this decomposition is a FRH decomposition of \(\phi_t\) such that:
\begin{enumerate}

\item each \(C(h_j)\) is one of types (d), (e) and (f) in Lemma \ref{lem2};
\item every boundary component of each \( \widetilde{C}(h_{j})\) is incompressible in \(\widetilde{C}(h_{j})\). 
\end{enumerate}

\end{lemm}

\begin{proo}
Suppose that \(W=  (\partial_{-}W \times I) \cup_{j=1}^{n}\widetilde{C}(h_{j})\) is a FRH decomposition of \(\phi_t\). By the definition of \(\phi_t\),  we have \(n \geq 2\) when \(\partial W = \varnothing\). Then the boundary of each \( \widetilde{C}(h_{j})\) is nonempty.
 Let \(l\) be the indexed link of \(\phi_t\), then
the intersection of \(l\) and each RH  associated to the decomposition \(W=  (\partial_{-}W \times I) \cup_{j=1}^{n}\widetilde{C}(h_{j})\) is a core of this RH.

\emph{Necessity}. Suppose that  \(\phi_t\) is related to an  incompressible torus decomposition
 \(W=  (\partial_{-}W \times I) \cup_{j=1}^{n}\widetilde{C}(h_{j})\).
Since \(l\) is knotted in \(W\) and \(l \cap h_j\) is a core of \(h_j\), each  \(C(h_j)\) cannot be of type (a), (b) or (c) by Lemma \ref{lem2}. If there is a  \(C(h_{j_{0}})\)  of type (g), then by the irreducibility of \(W\), \(W \cong \mathbb{RP}^3\). Moreover, it is easy to prove that each boundary component of \(C(h_{j_{0}})\) is compressible in \(W\). By the definition of incompressible torus decompositions, each boundary component of \(C(h_{j_{0}})\) must bound either a \(0\)-FRH or \(2\)-FRH, which contradicts that  \(n \geq 2\). Therefore, each \(C(h_j)\) is one of types (d), (e) and (f) in Lemma \ref{lem2}.
   Let \(T\) be a boundary component of \( \widetilde{C}(h_j)\), then \(T\) does not bound any \(0\)-FRH or \(2\)-FRH. By the definition of incompressible torus decompositions, \(T\) is incompressible in \(W\), then \(T\) is incompressible in \(\widetilde{C}(h_{j})\). 

\emph{Sufficiency}. Suppose that the FRH decomposition \(W=  (\partial_{-}W \times I) \cup_{j=1}^{n}\widetilde{C}(h_{j})\) of \(\phi_t\) satisfies conditions (1) and (2).
Let \(T'\) be 
 a boundary component  of a \(1\)-FRH in this decomposition that does not bound a 0-FRH or 2-FRH, then  \(T'\) is a boundary component of  \( \widetilde{C}(h_{j'})\) for some \(j' \in \{1, \cdots, n\}\). By the condition (2), \(T'\) is incompressible in \(\widetilde{C}(h_{j'})\). Since \(W=(\partial_{-}W \times I) \cup_{j=1}^{n}\widetilde{C}(h_{j})\) and every boundary component of each \( \widetilde{C}(h_{j})\) is incompressible in \(\widetilde{C}(h_{j})\), \(T'\) is incompressible in \(W\). 

Let \(l_j=l \cap \widetilde{C}(h_{j})\) for  \(j=1, \cdots, n\). Note that 
the intersection of \(l_j\) and each RH  associated to \(\widetilde{C}(h_{j})\) is a core of this RH.
According to the proof of Lemma \ref{lem3}, \(\widetilde{C}(h_{j})\) is a Seifert manifold, and \(l_j\) consists of some fibers of a Seifert fibering of \(\widetilde{C}(h_{j})\). Then \(l_j\) is knotted in \(\widetilde{C}(h_{j})\). Suppose that there is a knot \(k\) of \(l\) bounding a disk in \(W\), and \(k\) is a knot of some \(l_{j''}\) for \(j'' \in \{1, \cdots, n\}\). Since every boundary component of each \( \widetilde{C}(h_{j})\) is incompressible in \(\widetilde{C}(h_{j})\), \(k\)
must bound a disk in \(\widetilde{C}(h_{j''})\), which contradicts that \(l_{j''}\) is knotted in \(\widetilde{C}(h_{j''})\). Therefore, \(l\) is knotted in \(W\).
Then \(\phi_t\) is related to an incompressible torus decomposition
\(W=  (\partial_{-}W \times I) \cup_{j=1}^{n}\widetilde{C}(h_{j})\).

\end{proo}


\begin{prop}
\label{prop2}

Let \(\phi_t\) be an NMS flow on  a graph manifold \(W\)  related to  incompressible torus decompositions. 
Then there is a JSJ decomposition \(W= M_1 \cup \cdots \cup M_s\) with the JSJ tori set \(\mathcal{T}\) such that:

\begin{enumerate}
\item \(\phi_t\) is transverse to each 
  \(T \in \mathcal{T}\);

\item    \(\phi_t |_{M_{i}}\) is an NMS flow on \(M_i\) related to incompressible torus decompositions.
\end{enumerate}

\end{prop}

\begin{proo}
If \(W\) is a Seifert manifold, then the conclusion in this case is antomotially correct. From now on, we suppose that \(W\) is not a Seifert manifold, i.e. \( \mathcal{T} \neq \varnothing\).

Since \(\phi_t\) admits saddle closed oribits, any FRH decomposition of \(\phi_t\) must contain \(1\)-FRHs.
Suppose that \(\phi_t\) is related to an incompressible torus decomposition \( W= (\partial_{-}W \times I) \cup_{j=1}^{n} \widetilde{C}(h_{j})\). 
Since the indexed link of \(\phi_t\) is knotted in \(W\), every \(C(h_j)\) cannot be of type (a), (b) or (c). 
Since \(W\) is not a Seifert manifold, \(C(h_j)\) cannot be of type (g). Thus  each \(C(h_j)\) is one of types (d), (e) and (f).

If \(n=1\) and \(\partial W \neq \varnothing\), then  \( W= (\partial_{-}W \times I) \cup \widetilde{C}(h_{1}) \cong \widetilde{C}(h_{1})\).
By Lemma \ref{lem3} and the irreducibility of \(W\), \(W\)
is a Seifert manifold. If \(n=1\) and \(\partial W = \varnothing\),
then similar to Lemma \ref{lem3}, we can prove that  \(W\)  is homeomorphic to  
\(M(0,0;\frac{q_1}{p_1}, \frac{q_2}{p_2}, \frac{q_3}{p_3})\) or \(M(-1,0; \frac{q'_1}{p'_1}, \frac{q'_2}{p'_2} )\). Therefore, when \(n=1\), \(W\)
is a Seifert manifold, which contradicts to the supposition.
Then \(n \geq 2\),
which implies that the boundary of every \( \widetilde{C}(h_{j})\) is nonempty. By the definition of \(\phi_t\), the boundary components of each \( \widetilde{C}(h_{j})\) are incompressible in \(W\). By Lemma \ref{lem3},  each \( \widetilde{C}(h_{j})\) is a Seifert manifold. Then there is a JSJ decomposition \(W= M_1 \cup \cdots \cup M_s\)  with the JSJ tori set \(\mathcal{T}\)
such that for any 
 \(T \in \mathcal{T}\), there is \(k\in \{1, \cdots, n\}\) such that
\(T\) is a boundary component of   \(\widetilde{C}(h_{k})\). 
Then \(\phi_t \) is transverse to \(T\).
Moreover, the indexed link of \(\phi_t|_{M_i}\) 
is knotted  in \(M_i\) and must contain index-\(1\) knots.
It is easy to observe that \(\phi_t |_{M_{i}}\) is an NMS flow related to incompressible torus decompositions.
Proposition \ref{prop2} is proved.

\end{proo}

\section{Indexed links of NMS flows on ordinary graph manifolds} \label{s.indexed links}

 In Section \ref{s.int}, we introduced the concept of ordinary graph manifolds. Recall that  a closed graph manifold \(W\) is called an \emph{ordinary graph manifold} if:
\begin{enumerate}
\item each Seifert piece \(M_i\) of \(W\)  admits a unique Seifert fibering up to isotopy;
\item the base orbifold of \(M_i\)  is orientable, and \(M_i\) does not admit any singular fiber with slope-\(\frac{q}{2}\) where \(q\) is coprime to \(2\);
\item \(W\) is not homeomorphic to \(M(0,0; \frac{q_1}{p_1}, \frac{q_2}{p_2}, \frac{q_3}{p_3})\) (\(p_1 , p_2, p_3 > 1\)).
\end{enumerate}

If \(W\) is an ordinary graph manifold, then we claim that \(W\) cannot admit  a Seifert fibering whose base orbifold is a $2$-sphere with $k$ singularities ($k\leq 3$).
 Now, we will prove this claim.
We assume by contradiction  that \(W\) is  homeomorphic to some \(M(0,0; \frac{q_1}{p_1}, \frac{q_2}{p_2}, \frac{q_3}{p_3})\).
Without loss of generality, we assume that \(p_1, p_2, p_3 \geq 1\). According to the condition (3) in the definition of ordinary graph manifolds, the case that \(p_1 , p_2, p_3>1\) is impossible. For other cases, \(W\) contains at most two singular fibers, then \(W\) admits some different Seifert fiberings up to isomorphism (see Theorem 2.3 in Hatcher \cite{Ha}), which contradicts the first condition in the definition of ordinary graph manifolds. Therefore,  \(W\) cannot be homeomorphic to  \(M(0,0; \frac{q_1}{p_1}, \frac{q_2}{p_2}, \frac{q_3}{p_3})\).

\par
From now on, we discuss the relationship between NMS flows and NMS flows related to incompressible torus decompositions. 
\begin{theorem}\label{t.genindlink1}
Let $l$ be an indexed link in an ordinary graph manifold $W$. Then $l$ is the indexed link of
an NMS flow \(\phi_t\)  on $W$ if and only if
there is an NMS flow $\psi_t$ on $W$ related to incompressible torus decompositions such that
$l$
can be obtained from the indexed link $l'$ of  $\psi_t$  by  applying finitely many steps of  operations in Operation A.

\end{theorem}

The proof of Theorem \ref{t.genindlink1} is very complicated. 
For the convenience of the reader, before the complete proof, we outline the idea of the proof.

\textbf{Proof ideas:}

 \textbf{Necessity}. Let \(l\) be the indexed link
of an NMS flow \(\phi_t\)  on $W$.
\begin{enumerate}
\item By Theorem \ref{t.thm1},
we can get a FRH decomposition of \(\phi_t\), which does not contain any \(1\)-FRH of type (g) in Lemma \ref{lem2}.
 \item Note that a \(1\)-FRH of type (a), (b) or (c) is a connected sum of two \(3\)-manifolds.
By the irreducibility of \(W\), we can get a
 FRH decomposition of  some NMS flow \(\phi_t^0\) on \(W\),
which does not contain   any \(1\)-FRH of type (a), (b), (c) or (g). Moreover,  \(l\) is obtained from the indexed link \(l_0\) of \(\phi_t^0\) by applying finitely many steps of Operations  \uppercase\expandafter{\romannumeral1}-\uppercase\expandafter{\romannumeral4} in Operation A.

\item By  replacing some solid tori with \(0\)-FRHs or  \(2\)-FRHs, 
 we can construct an NMS flow \(\psi_t\) on \(W\) related to incompressible torus decompositions
  such that \(l_0\) is obtained from the indexed link \(l'\) of \(\psi_t\) by applying finitely many steps of Operations  \uppercase\expandafter{\romannumeral5}-\uppercase\expandafter{\romannumeral7} in Operation A.
\end{enumerate}

\textbf{Sufficiency}. 
Let \(l_1\) be the indexed link of an NMS flow \(\phi_t^1\) on \(W\), and \(l_0\) be an indexed link in \(W\).

\begin{enumerate}
\item If \(l_0\)  is obtained from \(l_1\) by applying one of Operations  \uppercase\expandafter{\romannumeral1}-\uppercase\expandafter{\romannumeral5} in Operation A, then by doing  connected sum \(W \sharp S^3\), we can construct an NMS flow   on \(W\) with indexed link \(l_0\).

\item  If \(l_0\)  is obtained from \(l_1\) by applying Operation \uppercase\expandafter{\romannumeral6} or \uppercase\expandafter{\romannumeral7} in Operation A, then 
by changing the flow on a filtrating neighborhood of an attracting or repelling closed orbit of  \(\phi_t^1\), we can construct a new NMS flow   on \(W\) with indexed link \(l_0\).

\end{enumerate}

\begin{proo}[Proof of Theorem \ref{t.genindlink1}]
\emph{Necessity}.
Let \(l\) be the indexed link
of an NMS flow \(\phi_t\)  on $W$.
By  Theorem \ref{t.thm1},  we can obtain a FRH decomposition \(W=  \cup_{i} C_{i}\)  of \(\phi_t\), then
the intersection of \(l\) and each RH  associated to the decomposition \(W=  \cup_{i} C_{i}\) is a core of this RH.

If the decomposition   \(W=  \cup_{i} C_{i}\) does not contain any \(1\)-FRH,
then \(W\) is homeomorphic to one of \(T^2 \times I\), \(S^1 \times D^2\), \(S^3\) and a lens space, which contradicts to the definition of  ordinary graph manifolds.
Thus, this decomposition must contain \(1\)-FRHs.
 In addition, since \(W\) cannot be homeomorphic to \(\mathbb{RP}^3\), this decomposition cannot contain any \(1\)-FRH of type (g) in Lemma \ref{lem2}.  
Let \(C_{i_{0}}\) be a \(1\)-FRH associated to a \(1\)-RH \(h\) in this FRH decomposition.

\par Suppose that \(C_{i_{0}}\cong (T_{1}\times I) \sharp (T_{2}\times I)\), where \(T_{1}\) and \(T_{2}\) are two tori. Then a \(2\)-sphere \(S^{2}\) associated to the connected sum separates both of \(C_{i_{0}}\) and \(W\). Suppose that \(C_{i_{0}}= A_{1}\cup_{S^{2}}B_{1}\) and \(W= A\cup_{S^{2}}B\), where \(A_{1} \cup_{S^{2}} B^{3} \cong T_{1}\times I\), \(B_{1} \cup_{S^{2}} B^{3} \cong T_{2}\times I\), \( A_{1}\subseteq A\), and  \(B_{1}\subseteq B\). Then 

\begin{equation}
\begin{split}
\centering
W &= A\cup_{S^{2}}B =(A\smallsetminus {\rm Int}A_{1})\cup_{\partial A_{1}\smallsetminus S^2}A_{1}\cup_{S^{2}}B_{1}\cup_{\partial B_{1}\smallsetminus S^2}(B\smallsetminus {\rm Int}B_{1})\\
&=(A\smallsetminus {\rm Int}A_{1})\cup_{T_{1}\times \partial I} ((T_{1} \times I) \sharp (T_{2} \times I)) \cup_{T_{2}\times \partial I}(B\smallsetminus {\rm Int}B_{1})\nonumber.
\end{split}
\end{equation}
By the irreducibility of \(W\), either \((A\smallsetminus {\rm Int}A_{1})\cup_{T_{1}\times \partial I} (T_{1} \times I) \cong S^{3}\) or \( (T_{2} \times I) \cup_{T_{2}\times \partial I}(B\smallsetminus {\rm Int}B_{1}) \cong S^{3}\). Without loss of generality, we assume that \( (T_{2} \times I) \cup_{T_{2}\times \partial I}(B\smallsetminus {\rm Int}B_{1}) \cong S^{3}\), then \((A\smallsetminus {\rm Int}A_{1})\cup_{T_{1}\times \partial I} (T_{1} \times I) \cong W\). 
Since both \(A\smallsetminus {\rm Int}A_{1}\) and \(B\smallsetminus {\rm Int}B_{1}\) admit FRH decompositions induced by \(W= \cup_{i} C_{i}\),  we get  a new FRH decomposition of \(W\) which contains fewer \(1\)-FRHs of type (a), and get  a FRH decomposition of \(S^{3}\). By  Theorem \ref{t.thm1}, we can construct  an NMS flow on \(W\) (resp. \(S^3\)) such that the above FRH decomposition of \(W\) (resp.  \({S}^3\)) is a FRH decomposition of this flow. Let \(l'_1\) (resp. \(l'_2\)) be the indexed link of the above NMS flow on \(W\) (resp. \(S^3\)). It is easy to observe that \(l=l'_{1}\cdot l'_{2}\cdot u\), where \(u\) is a core of \(h\). Namely,  \(l\) is obtained from \(l'_{1}\) and \(l'_{2}\) by applying Operation \uppercase\expandafter{\romannumeral1} in Operation A.

For  this new FRH decomposition of \(W\), we 
continue to
discuss its \(1\)-FRHs of type (a) in the same way. After finitely many steps, we end up with a  FRH decomposition of an NMS flow on \(W\), which does not contain any \(1\)-FRH of type (a).  Similarly, it is easy to observe that  discussing the \(1\)-FRH of type (b) corresponds to  Operation \uppercase\expandafter{\romannumeral2} or \uppercase\expandafter{\romannumeral3}, and discussing the \(1\)-FRH of type (c)  corresponds to  Operations \uppercase\expandafter{\romannumeral4}. 
Therefore, 
 we  can finally obtain a FRH decomposition \(W= \cup_{w} C_{w}^{0}\) of some NMS flow  \(\phi_t^0\) on \(W\), which does not contain any  \(1\)-FRH of type (a), (b), (c) or (g). Let  \(l_0\) be the indexed link of \(\phi_t^0\), then \(l\) is obtained from \(l_0\) by applying finitely many steps of Operations  \uppercase\expandafter{\romannumeral1}-\uppercase\expandafter{\romannumeral4} in Operation A.
 
\begin{clai}
\label{claim1}
There is an NMS flow \(\psi_t\) on \(W\) related to incompressible torus decompositions such that
\(l_0\) is obtained from the indexed link \(l'\) of  \(\psi_t\)  by applying finitely many steps of Operations \uppercase\expandafter{\romannumeral5}-\uppercase\expandafter{\romannumeral7} in Operation A.
\end{clai}

\begin{proo}
If the  decomposition  \(W= \cup_{w} C_{w}^{0}\)  does not contain any \(1\)-FRH,
then \(W\) is homeomorphic to  one of \(T^2 \times I\), \(S^1 \times D^2\), \(S^3\) and a lens space. 
If 
 the decomposition \(W= \cup_{w} C_{w}^{0}\)   contains only one \(1\)-FRH,  then similar to Lemma \ref{lem3}, we can prove that 
\(W\) is homeomorphic to  
\(M(0,0;\frac{q_1}{p_1}, \frac{q_2}{p_2}, \frac{q_3}{p_3})\) or \(M(-1,0; \frac{q'_1}{p'_1}, \frac{q'_2}{p'_2} )\). 
The above two cases 
contradict to the definition of \(W\). Therefore the decomposition \(W= \cup_{w} C_{w}^{0}\)   contains at least two \(1\)-FRHs.

Let  \(W= \cup_{w} C_{w}^{0} = \cup_{j=1}^{n} \widetilde{C}(h_{j})\) (\(n \geq 2\)), where \(C(h_{j})\) is a  \(1\)-FRH of type (d), (e) or (f) associated to a \(1\)-RH \(h_j\). 
Then
 the boundary of each \( \widetilde{C}(h_{j})\) is nonempty.
We say that a boundary component $T$ of some \( \widetilde{C}(h_{j'})\) is  a \emph{compressible torus of the FRH decomposition $W= \cup_{w} C_{w}^{0}$} if $T$  is compressible in \(W\).

%
%

\par Let \(y\) be the number of the compressible tori of the FRH decomposition \(W= \cup_{w} C_{w}^{0}\). If \(y =0\), then every boundary component of \(\widetilde{C}(h_{j})\) is incompressible in \(\widetilde{C}(h_{j})\) for every \(j=1, \cdots, n\). By Lemma \ref{lemITD}, 
  \(\phi_t^0\) is related to an incompressible torus decomposition \(W= \cup_{w} C_{w}^{0}\).
If \(y>0\), then there is some \( j_{0}\in \left \{ 1,\cdots, n \right \} \) such that \(\widetilde{C}(h_{j_{0}})\) supports compressible boundary components  in \(\widetilde{C}(h_{j_{0}})\). 
Since \(W\) is not homeomorphic to any lens space, by Lemma \ref{lem3},
 \(\widetilde{C}(h_{j_{0}})\) is homeomorphic to \( V_{1} \sharp V_{2}\) or \(V_1\), where \(V_1\) and \(V_2\) are two solid tori. Let \(m= l_0 \cap h_{j_0}\), then \(m\) is a core of \(h_{j_{0}}\) with index \(1\).

\par \textbf{Case 1.} Suppose that \(\widetilde{C}(h_{j_{0}}) \cong V_{1} \sharp V_{2}\).

Similar to the discussion of the \(1\)-FRH of type (a), we can get a  FRH decomposition of an NMS flow on \(W\) such that  the number of compressible tori of this FRH decomposition is less than \(y\). In addition, we get  a FRH decomposition of an NMS flow on \(S^{3}\). In fact, one of \(V_1\) and \(V_2\) is a FRH in this decomposition of \(W\), and the other is a FRH in this  decomposition of \(S^3\). Without loss of generality, we assume that \(V_1\) is a FRH in this  decomposition of \(W\).
Let \({l}_1\) and \({l}_2\) be the indexed links of the above NMS flows on \(W\) and \({S}^3\) respectively. Let \(k_1 = l _1 \cap V_1\) and \(k_2 =l_2 \cap V_2\), then \(k_1\) is a core of \(V_1\) and \(k_2\) is a core of \(V_2\). Moreover, each of the indices of \(k_1\) and \(k_2\) is either \(0\) or \(2\).

Due to the proof of Lemma \ref{lem3},  \(C(h_{j_{0}})\cong F \times S^1\), and \(m= \{\ast\} \times S^1\), where \(F\) is a pair-of-pants and \(\ast \in {\rm Int} F\). Moreover,
\(\widetilde{C}(h_{j_{0}})= C(h_{j_{0}}) \cup R_{0}\), where \(R_0\) is a \(0\)-FRH or \(2\)-FRH. Let \(\partial_{0}\), \(\partial_{1}\) and \(\partial_{2}\) be the boundary components of \(F\). Suppose that \(\partial_{0} \subset \partial R_0\), \( \partial _{1} \subset \partial V_1 \) and \( \partial _{2} \subset \partial V_2 \).
Then  for each \(p\in \partial_{0}\), \(\{p\} \times S^1\) bounds a meridian disk of \(R_0\).
 Let \(c(R_0)= l_0 \cap R_0\), then it is a core of \(R_0\).

Let \(\alpha\) be an essential arc in \(F\) such that its end points \(p_1\) and \(p_2\) are contained in 
\(\partial_{0}\) (see Figure \ref{Fgthm2_1} (a)), then \(\{p_1\} \times S^1\)  (resp. \(\{p_2\} \times S^1\)) bounds a  meridian disk \(D_1\) (resp. \(D_2\)) of \(R_0\). We get \(S^2 \cong D_1 \cup (\alpha \times S^1) \cup D_2\) which is a separating \(2\)-sphere  associated to the connected sum of \(V_1 \sharp V_2\). It is easy to observe that \(c(R_{0})\simeq k_1 \sharp k_2\) whose index is equal to either \({\rm Ind}(k_1)\) or  \({\rm Ind}(k_2)\), and \(m\) is a meridian of \(k_1 \sharp k_2\) (see Figure \ref{Fgthm2_1} (b)). 
Thus, \(l_0\) is obtained from \(l_{1}\) and \(l_{2}\) by applying Operation \uppercase\expandafter{\romannumeral5} in Operation A. 

\begin{figure}[htbp]
\centering
\subfigure[ ]{\includegraphics[width=0.38\textwidth]{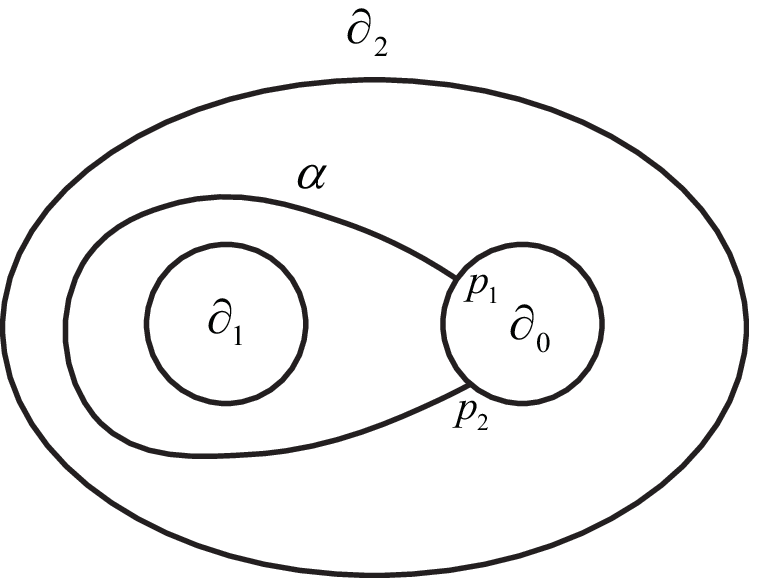}}
\hspace{.60in}
\subfigure[ ]{\includegraphics[width=0.45\textwidth]{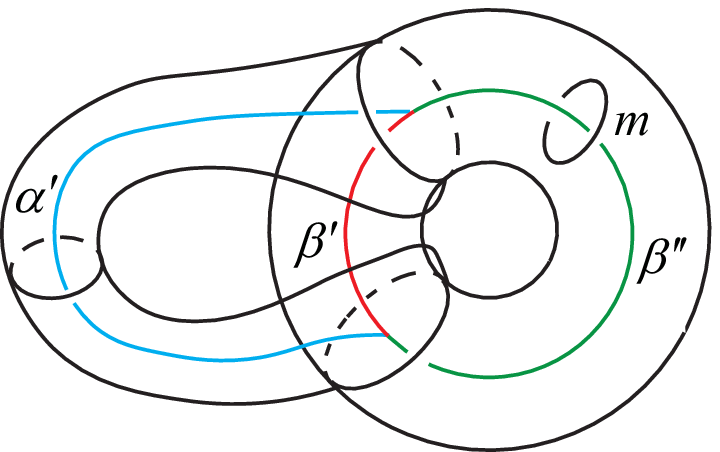}}
\hspace{.60in}

\caption{\(\beta' \cup \beta'' \simeq c(R_0)\), \(\beta' \cup \alpha' \simeq k_1\), and \(\beta'' \cup \alpha' \simeq k_2\). }
\label{Fgthm2_1}
\end{figure}

\par \textbf{Case 2.} Suppose that \(\widetilde{C}(h_{j_{0}}) \cong V_1\). 

\par According to \(\partial _{-} (\widetilde{C}(h_{j_{0}})) = \varnothing \) or \(T^{2}\), we replace \(\widetilde{C}(h_{j_{0}})\) with a \(0\)-FRH or \(2\)-FRH \(R\). Then we get a  FRH decomposition of an NMS flow on \(W\), such that the number of compressible tori of this  FRH decomposition is less than \(y\). 
Let \(l_{1}\) be the indexed link of the above NMS flow, and \(k= l_1 \cap R\), then \(k\) is  a core  of \(R\) with index \(2\) or \(0\).

\par Suppose that  \(C(h_{j_{0}})\) is of type (d), then \(\widetilde{C}(h_{j_{0}}) = R_1 \cup C(h_{j_0}) \cup R_2\), where \(R_1\) and  \(R_2\) are FRHs which are not \(1\)-FRHs.  Let \(N'= R_1 \cup C(h_{j_0})\) and \(N''=C(h_{j_0}) \cup R_2\).
By the proof of Lemma \ref{lem3},  one of the following holds:
\begin{itemize}

\item One of \(N'\) and \(N''\) is the same to  \(\widetilde{C}(h)\) of type (5) in Lemma \ref{lem3}, and the other is not.
\item \(N'\) or \(N''\) is homeomorphic to \(T^2 \times I\).
\end{itemize}

For the first case, similar to Case 1, we can get an indexed link \(l_{2}\) of an NMS flow on \(S^{3}\) such that  \(l_0\) is obtained from \(l_{1}\) and \(l_{2}\) by applying Operation \uppercase\expandafter{\romannumeral5}. For the second case, we may assume that \(N' \cong T^2 \times I\).
Then  \(\widetilde{C}(h_{j_{0}})\) contains three indexed knots of \(l_0\): \(k_1\), \(k_2\) and \(m\), where \(k_1\) (resp. $k_2$) is a core of \(R_1\) (resp. $R_2$). 
Moreover, \(k_{1}\) and \(m\) are two parallel \(\left ( p,q \right )\)-cables of \(k_2\), and \(k_2\) is  a core of the solid torus \(\widetilde{C}(h_{j_{0}})\). The indices of \(k_1\) and \(k_{2}\) are either \(0\) or \(2\), and one of them is equal to \({\rm Ind}(k)\).  Namely, \(l_0\) is obtained from \(l_1\) by applying Operation \uppercase\expandafter{\romannumeral6}.

\par Suppose that \(C(h_{j_{0}})\) is of type (e), then we can get \(F \times S^1\) by removing a small open tubular neighborhood  of \(m\) from \(C(h_{j_{0}})\), where \(F\) is a pair-of-pants. 
Similar to the above discussion, it is easy to observe that \(l_0\) is obtained from \(l_{1}\) by applying Operation \uppercase\expandafter{\romannumeral7}.

In both Case 1 and Case 2, we get a FRH decomposition of an NMS flow on \(W\) 
such that  the number of compressible tori of this FRH decomposition is less than \(y\). For this FRH decomposition, we continue to discuss its
compressible tori in the same way. After finitely many steps, we end up with a FRH decomposition  
 of some NMS flow \(\psi_t\) on \(W\) such that there is no
  compressible torus of this decomposition. 
 By Lemma \ref{lemITD},
\(\psi_t\) is an NMS flow related to   incompressible torus decompositions. Let \(l'\) be the indexed link of \(\psi_t\), then  \(l_0\) is obtained from  \(l'\)  by applying finitely many steps of Operations \uppercase\expandafter{\romannumeral5}-\uppercase\expandafter{\romannumeral7} in Operation A. Claim  \ref{claim1} is proved.

\end{proo}
\par By Claim \ref{claim1}, the necessity of Theorem \ref{t.genindlink1} is proved. 
\\
\par \emph{Sufficiency}.
Let \(l_1\) and \(l_2\) be two indexed links of   two NMS flows on \(W\) and \(S^3\) respectively. Let \(W=  \cup_{j=1}^{s}C_{j}\) and \({S}^3=\cup_{j=1}^{w}C'_{j}\) be two  FRH decompositions of the above flows respectively.  Let \(T_{1}\) (resp. \(T_2\)) be a boundary component of \(C_1\) (resp. \(C'_1\)), and let \(l_0\) be an indexed link in \(W\).

If \(l_0\) is obtained from \({l}_1\) and \({l}_2\) by applying Operation \uppercase\expandafter{\romannumeral1} in Operation A, we repalce \({C}_1\) by \({C}_1 \cup (T_1 \times I)\) and repalce \(C'_{1}\) by \(C'_{1} \cup (T_2 \times I)\). Then 
\begin{center} 
\(W \cong W \sharp {S}^3 \cong  \left( T_1 \times I \right) \sharp \left( T_{2} \times I \right) \cup \left(\cup_{j=1}^{s}C_{j}\right) \cup \left(\cup_{j=1}^{w}C'_{j}\right)\).
\end{center}
By considering \( (T_1 \times I) \sharp (T_2 \times I)\) as a \(1\)-FRH of type (a) in Lemma \ref{lem2},
we get a new FRH decomposition of \(W\).
By
Theorem \ref{t.thm1}, we can construct  an NMS flow on \(W\) whose  indexed link is \(l_0\).

Similarly,  if  \(l_0\) is obtained from \({l}_1\) and \({l}_2\) by applying Operation
 \uppercase\expandafter{\romannumeral2}\ or \uppercase\expandafter{\romannumeral3}, then by using the \(1\)-FRH of type (b), we can construct  an NMS flow on \(W\) whose  indexed link is \(l_0\).
 If  \(l_0\) is obtained from \({l}_1\) and \({l}_2\) by applying  Operation
 \uppercase\expandafter{\romannumeral4}, then by using the \(1\)-FRH of type (c), we can construct an NMS flow on \(W\) whose  indexed link is \(l_0\). In addition, if  \(l_0\) is obtained from $l_1$ and $l_2$ by applying Operation
 \uppercase\expandafter{\romannumeral5}, then by using \(\widetilde{C}(h)\) of type (5) in Lemma \ref{lem3}, we can construct  an NMS flow on \(W\) whose  indexed link is \(l_0\).
If  \(l_0\) is obtained from $l_1$ by applying Operation
 \uppercase\expandafter{\romannumeral6}\ or \uppercase\expandafter{\romannumeral7}, then by using \(\widetilde{C}(h)\) of type (4) in Lemma \ref{lem3}, we can construct  an NMS flow on \(W\) whose  indexed link is \(l_0\).

Let $\psi_t$ be an NMS flow on $W$ related to incompressible torus decompositions, and \(l'\) be the indexed link of \(\psi_t\). Suppose that 
 an indexed link $l$ in \(W\)
can be obtained from  $l'$ by applying  finitely many steps of  operations in Operation A.
According to the above discussion, there is an NMS flow \(\phi_t\) on \(W\) such that \(l\) is the indexed link of \(\phi_t\). Sufficiency is proved, and therefore the proof of Theorem \ref{t.genindlink1} is complete.

\end{proo}

Let \(l_1\) be the indexed link of an NMS flow \(\phi_t\) on an ordinary graph manifold \(W\). Let
 \(N \cong T^2 \times I\) be a submanifold of \(W\) such that \(\phi_{t} |_{N}\) is topologically equivalent to the flow induced by \(\frac{\partial}{\partial x}\) 
along the \(I\) direction. Let \(V\) be a filtrating neighborhood of an attracting or repelling closed orbit of \(\phi_t\).
According to the proof of Theorem \ref{t.genindlink1}, Operations \uppercase\expandafter{\romannumeral1}\ and  \uppercase\expandafter{\romannumeral2} are corresponding to change the NMS flow \(\phi_{t} |_{N}\), and Operations \uppercase\expandafter{\romannumeral3}-\uppercase\expandafter{\romannumeral7} are corresponding to change the NMS flow \(\phi_{t} |_{V}\).

\section{Lyapunov graphs}\label{ss.Lya}

\begin{defi}
A \emph{generalized graph} $G=(V, V', E)$ is a topological space obtained by connecting a set $V\sqcup V'$ of points and a set $E$ of edges, such that 
\begin{enumerate}
\item
 each point in $V'$ is connected by only one edge;
\item 
 $(V\sqcup V', E)$ is a \emph{graph}, i.e., each edge in $E$ is connected to two distinct points in $V\sqcup V'$.
\end{enumerate}
\end{defi}
For a  generalized graph $G=(V, V', E)$,
we call the points in $V$ the \emph{vertices} of $G$, call
the points in $V'$ the \emph{ends} of $G$. For a vertex $v$ of  $G$, the \emph{degree} of $v$ denotes  the number of the edges connecting $v$. 
$G$ is called a \emph{tree} if its first Betti number \(\beta _{1}(G)\) is equal to \(0\).
\begin{rema}
 
 The generalized graph is actually a graph, except that the ``vertices" of the graph are divided into two categories: the vertices of the generalized graph and the ends of the generalized graph. Moreover, when we consider an end of the generalized graph as a vertex of a graph (in standard concept), its degree is $1$.
 

\end{rema}

Let $G_i$  be a generalized graph, and $p_i$ be an end of $G_i$ connected by the edge $e_i$ in $G_i$ for $i=1,2$. 
In this paper, ``\emph{gluing $G_1$ and $G_2$ along $p_1, p_2$}'' means that glue $G_1$ and $G_2$ together by gluing $p_1$ and $p_2$, and $e_1, e_2, p_1, p_2$ form an edge in the resulting generalized graph $G$. Moreover, we  specify that the vertices of $G$ consist of the vertices of $G_1 \sqcup G_2$. In addition, 
if a generalized graph $G'_1$ is obtained by cutting $G_1$ along a point in an edge, then we  specify that the vertices of $G'_1$ consist of the vertices of $G_1$.

The Lyapunov graph was first  used by  Franks  \cite{Fr} to classify nonsingular Smale flows  on $S^3$. Let $\phi_t$ be a smooth flow on  a compact manifold $M$ with a Lyapunov function \(f: M\rightarrow \mathbb{R}\), where $f$ maps each  component of $\partial M$ to a constant. A \emph{Lyapunov graph} is an oriented generalized graph by identifying each connected component of \(f^{-1}(c)\) to a point for each \(c \in \mathbb{R}\), where the components of the level sets of $f$  that contain closed orbits produce the vertices of $L$, and the boundary components of $M$ produce the ends of $L$.
Moreover,  each edge is oriented by the flow direction. 
In fact, we can also define the Lyapunov graph in an abstractive way (see Franks \cite{Fr}):

\begin{defi}
An \emph{abstract Lyapunov graph} is a finite, connected, oriented generalized  graph \(L\) which satisfies the following conditions:
\begin{enumerate}

\item \(L\) possesses no oriented cycles;

\item each vertex of \(L\) is labeled with a chain recurrent flow on a compact space.
\end{enumerate}
\end{defi}

Let \(L\) be an abstract Lyapunov graph, and \(v\) be a vertex of  \(L\). We denote by  \(e_{v}^{-}\) (resp.\(e_{v}^{+}\))  the number of incoming (resp. outgoing) edges connecting \(v\). If \(e_{v}^{-}\cdot e_{v}^{+} \neq 0\), then  we call \(v\) a \emph{saddle vertex}. If \(e_{v}^{-}=0\) (resp. \(e_{v}^{+} =0\)), we call \(v\) a \emph{source (resp. sink) vertex}. 

\par
\begin{figure}[htbp]
\centering
\includegraphics[scale=0.6]{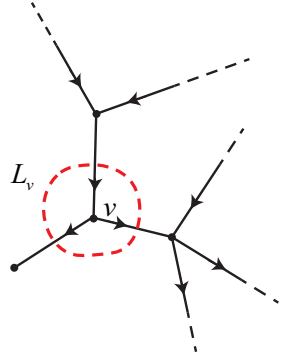}
\caption{\(L_{v}\) is the star neighborhood of \(v\).}
\label{FgLyaneg}
\end{figure}

 By cutting \(L\) along the midpoint of each edge connecting two vertices, we obtain some connected generalized graphs, each of which contains only one vertex. 
We call them \emph{the star neighborhoods of vertices} (see Figure \ref{FgLyaneg}). In particular, if \(L\) contains only one vertex \(v\), then we say that \(L\) is the star neighborhood of \(v\).

\begin{prop}\label{p.ends}
Let \(L\) be a Lyapunov graph of an NMS flow \(\phi_t\) on an orientable \(3\)-manifold, then   each of the  source and sink vertices is a degree $1$ vertex.
\end{prop}
\begin{proo}
Let 
\(v\) be a sink (resp. source) vertex of \(L\).
Note that the star neighborhood of \(v\) in \(L\) corresponds to a filtrating neighborhood of an attracting (resp. repelling) closed orbit of \(\phi_t\).
Since any filtrating neighborhood of an attracting (resp. repelling) closed orbit of \(\phi_t\) is a solid torus, \(v\) connects only one edge in \(L\). 

\end{proo}

\begin{figure}[htbp]
\centering
\includegraphics[scale=0.6]{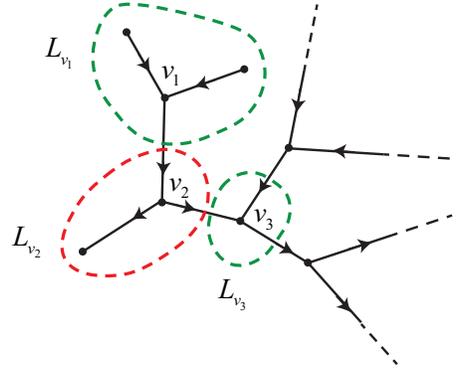}
\caption{\(L_{v_{i}}\) is the block associated to  \(v_i\) for \(i=1,2,3\).}
\label{FgLysub}
\end{figure}

 Let \(L\) be an abstract Lyapunov graph with saddle vertices such that  each of the source and sink vertices is a degree $1$ vertex.
 By cutting \(L\) along the midpoints of the edges  connecting two saddle vertices, we obtain some generalized graphs, each of which contains one saddle vertex. We call them  \emph{the blocks associated to saddle vertices} (see Figure \ref{FgLysub}). In particular, if \(L\) contains only one saddle vertex \(v\), then we say that \(L\) is the block associated to \(v\).

Let \(L'\) and \(L''\) be two finite oriented (possibly disconnected) generalized  graphs. If \(L''\) is attached to \(L'\) by an injective  map \(\varphi: \partial_{-}L'' \rightarrow \partial_{+}L'\), then we write the resulting generalized graph as \(L'+L''\), where \(\partial_{-}L''\) denotes the union of incoming ends of \(L''\),  and \(\partial_{+}L'\) denotes the union of outgoing ends of \(L'\). 
 
\begin{lemm}
\label{lemorder}
Let \(L\) be an abstract Lyapunov graph with saddle vertices such that  each of the source and sink vertices is a degree $1$ vertex.
Then \(L\) admits a decomposition \(L=(\partial_{-}L \times I)
+L_{1} +  \cdots + L_{n}\), where  \(L_1, \cdots, L_n\) are the blocks associated to saddle vertices of \(L\).

\end{lemm}

\begin{figure}[htbp]
\centering
\includegraphics[scale=0.6]{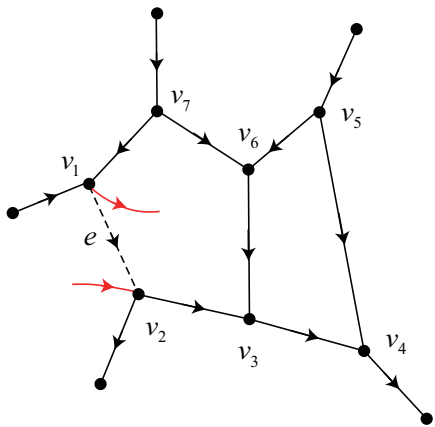}
\caption{Let \(L_{v_{i}}\) be the block associated to a saddle vertex \(v_{i}\). It is easy to observe that \(L' = (\partial _{-} L'\times I)+L_{v_{2}}+L_{v_{7}}+L_{v_{5}}+L_{v_{6}}+L_{v_{3}}+L_{v_{4}}+ L_{v_{1}}\), where \(v_{3}\) and \(v_{4}\) can be reached by \(v_{2}\) through oriented paths. Obviously, \(L=L_{v_{7}}+L_{v_{5}}+L_{v_{6}}+ L_{v_{1}}+L_{v_{2}}+L_{v_{3}}+L_{v_{4}}\). }
\label{FgdevomLy}
\end{figure}

\begin{proo}
Let \(w= \beta_{1}(L)\).
 By induction on the number of the saddle vertices, we can prove that \(L\)
admits  a decomposition when \(w=0\).
  Assume that \(w \geq 1\), and that an abstract Lyapunov graph with saddle vertices
admits  a decomposition if   each of the source and sink vertices is a degree $1$ vertex and its first Betti number is equal to \(w-1\).

By cutting \(L\) along the midpoint of  an oriented edge \(e\) in a circle  of \(L\), we get a generalized graph \(L'\) such that \(\beta_{1}(L')=w-1\) (see Figure \ref{FgdevomLy}). By the assumption of induction,
\(L'\) admits a decomposition 
 \(L' = (\partial_{-}L' \times I)+L_{1}+  \cdots + L_{n} \) where  \(L_1, \cdots, L_n\) are the blocks  associated to saddle vertices of \(L'\). In fact,  \(L_1, \cdots, L_n\) are also the blocks associated to saddle vertices of \(L\).

Let \(v_1\) and \(v_2\) be two saddle vertices of \(L\) such that
\(v_{1}\) reaches \(v_{2}\) through the oriented edge \(e\) in \(L\). Since \(L\) does not contain oriented circles, \(v_{2}\) cannot reach \(v_{1}\) through an oriented path in \(L'\). 
Suppose that \(L_i\) and \(L_j\) are the blocks associated to \(v_1\) and \(v_2\) respectively.
  If \(j>i\), then \(L\) admits a decomposition \(L= (\partial_{-}L \times I)+ L_{1}+  \cdots + L_{n}\). Otherwise, let \(S_{v_2}\) be the set of
the blocks  associated to saddle vertices that can be reached by \(v_2\) through oriented paths in \(L'\). Of course, $L_j$ is in $S_{v_2}$.
We move   the blocks in \(S_{v_2}\)  to the end of the decomposition of $L'$  perservng their order, and replace \(\partial _{-} L'\times I\) with \(\partial _{-} L\times I\). Then we  get a decomposition of \(L\). Lemma \ref{lemorder}  is proved.

\end{proo}

\begin{defi}\label{d.graphori}
Let \(\mathcal{S}=\{L\}\) be 
the collection of all connected graphs that satisfy the listed properties:

\begin{enumerate}
\item  the degree of each vertex of \(L\) is $1$ or $3$;
\item  the edges of \(L\) connecting degree $1$ vertices are oriented, and \(L\) contains both source degree $1$ vertices and  sink degree $1$ vertices;
\item if  \(p\) is a separating point of \(L\), then each component of \(L\smallsetminus \{p\}\)
 must contain degree $1$ vertices.

\end{enumerate}
\end{defi}

\begin{lemm}
\label{orientation}
Let \(L\) be a graph in \(\mathcal{S}\). Then other edges of \(L\) can be oriented such that the resulting graph is  an abstract Lyapunov graph, where each of the source and   sink vertices is a degree $1$ vertex. 
\end{lemm}
\begin{proo}
Let $x\geq 1$ be the number of the edges of $L$.  It is easy to observe that
the conclusion of the lemma is correct when \(x=1\).
Now we assume that \(x \geq 2\), and that the graphs in \(\mathcal{S}\) with less than \(x\) edges
satisfy the conclusion of the lemma.
We choose a path \(\gamma\)  with no self intersection
from a source  degree $1$ vertex \(v\) to a sink degree $1$ vertex \(v'\), and endow \(\gamma\) with an orientation from \(v\) to \(v'\). 
By deleting \(v\), \(v'\) and all  edges of \(\gamma\), we get some connected graphs \(L_1, \cdots, L_m\). Let \(n\) be the number of the vertices of \(L_1 \cap \gamma\), and  \(v_1, \cdots, v_n\) be the vertices in \(L_1 \cap \gamma\). Let \(p\) be a separating point in \(L_1\).

Suppose that \(n=1\), then
the point in the edge \(e\) of \(L_1\) connecting \(v_1\) is separating in \(L\).
By the condition (3) in the definition of \(\mathcal{S}\), the degree $1$ vertices of \(L_1\) consist of \(v_1\) and at least one degree $1$ vertex of \(L\).
Thus we can always endow an orientation on \(e\) such that \(L_1\) contains both  source degree $1$ vertices and  sink degree $1$ vertices.
Since  \(p\) is a separating point in \(L_1\),  \(p\) is separating in \(L\). Suppose that \(L_1\smallsetminus \{p\} =L'_1 \sqcup L''_1 \), where \(v_1 \in L'_1\). Then   \(L''_1\) is a component of \(L\smallsetminus \{p\}\). Due to the definition of \(\mathcal{S}\), \(L''_1\) must contain degree $1$ vertices. Therefore, \(L_1 \in \mathcal{S}\).

Suppose that  \(n \geq 2\).
For every \(i, j\), we define \(v_i \prec v_j\) if \(v_i \) points to \(v_j\) in \(\gamma\). Obviously, \(\{v_1, \cdots, v_n\}\) is a totally ordered set associated to \(\prec\). We may assume that \(v_1 \prec \cdots \prec v_n\).
In \(L_1\), we endow an orientation on the edges connecting \(v_1, \cdots, v_n\) such that \(v_1\) is a source degree $1$ vertex of \(L_1\) and \(v_2, \cdots, v_n\) are \(n-1\) sink degree $1$ vertices of \(L_1\). 

If \(p\) is  separating in \(L\), then similar to the case that \(n=1\), we can prove that each component of \(L_1\smallsetminus \{p\}\) contains degree $1$ vertices. If  \(p\) is  non-separating  in \(L\), then there is a path \(\gamma'\)  in \(L_1\)  with no self intersection from \(v_{i'}\) to \(v_{j'}\)  for \(i',j' \in \{1,\cdots, n\}\) such that \(p \in \gamma'\). Then one component of \(L_1\smallsetminus \{p\}\) contains \(v_{i'}\), and the other component of \(L_1\smallsetminus \{p\}\) contains \(v_{j'}\). Namely, each component of \(L_1 \smallsetminus \{p\}\) contains degree $1$ vertices. Therefore, \(L_1 \in \mathcal{S}\).

Obviously,  the number of the edges in \(L_1\) is less than \(x\).
By  the inductive assumption,  other edges of \(L_1\) can be oriented such that the resulting graph is  an abstract Lyapunov graph,  where each of the source and   sink vertices is a degree $1$ vertex. 
Similar discussions 
also work for
 \(L_2, \cdots , L_m\), and therefore   \(L\) can be oriented such that the resulting graph is  an abstract Lyapunov graph, where each of the source and   sink vertices is a degree $1$ vertex. 

\end{proo}
 
\par Suppose that \(L\) is a Lyapunov graph of an NMS flow \(\phi_t\) on an orientable \(3\)-manifold \(M\). Then \(L\) naturally corresponds to a FRH decomposition of \(\phi_t\), as follows.
By a Lyapunov function associated to \(L\), the star neighborhood of each vertex in \(L\) corresponds to a filtrating neighborhood of a closed orbit.  It is easy to prove that every filtrating neighborhood 
  is a FRH.  
By Lemma \ref{lemorder}, we can get a gluing order of the filtrating neighborhoods of the closed orbits in \(\phi_t\), and then we construct
a FRH decomposition
 \(M=  (\partial_{-}M \times I) \cup_{i=1}^{w} C_{i}\) of \(\phi_t\). Conversely, according to the numbers of the connected components of \(\partial_{-}C_{i}\) and \(\partial_{+}C_{i}\), we obtain a Lyapunov graph of  \(\phi_t |_{C_i}\).  Based on  the gluing rules of the FRHs of the FRH decomposition \(M=  (\partial_{-}M \times I) \cup_{i=1}^{w} C_{i}\), we get exactly  the
 Lyapunov graph \(L\) of \(\phi_t\).

\begin{prop}
\label{ITD SM}
Let \(M\) be a Seifert piece of an ordinary graph manifold \(W\) with a genus \(g\) base orbifold, and \(\phi_t\) be an NMS flow on \(M\). Suppose that  \(\phi_t\) is related to an 
 incompressible torus decomposition
\((\partial_{-}M \times I) \cup_{j=1}^{n} \widetilde{C}(h_{j})\), then:

\begin{enumerate}

\item For every \(j\), \(C(h_j)\) is of type (d) in Lemma \ref{lem2},
 \(\partial (\widetilde{C}(h_{j})) \neq \varnothing\),
 and 
\(\widetilde {C}(h_{j})\) is a Seifert manifold with a unique Seifert fibering up to isotopy unless \(\widetilde {C}(h_{j}) \cong  T^2 \times I\).

\item \(M\) is obtained from \(\widetilde{C}(h_{1}), \cdots, \widetilde{C}(h_{n}) \)  by some gluing  homeomorphisms that preserve the corresponding regular Seifert fibers. 
\item Let \(L\) be the Lyapunov graph of \(\phi_t\) corresponding to this decomposition, then  \( \beta_1 (L)=g\).
\end{enumerate}

\end{prop}

\begin{proo}
Since \(M\) cannot be homeomorphic to \(\mathbb{RP}^3\),  the decomposition \((\partial_{-}M \times I) \cup_{j=1}^{n} \widetilde{C}(h_{j})\) cannot contain any \(1\)-FRH of type (g) in Lemma \ref{lem2}. Since the indexed link of \(\phi_t\) is knotted in \(M\), every \(C(h_j)\) cannot be of  type  (a), (b) or (c).

Suppose that \(n=1\), then \(M\cong \widetilde{C}(h_{1})\), so \(\widetilde {C}(h_{1})\) is a Seifert manifold with the unique Seifert fibering up to isotopy. If \(\partial M = \varnothing\), then \(W \cong M\). Similar to Lemma \ref{lem3}, we can prove that \(W\) is homeomorphic to \( M(0,0;\frac{q_1}{p_1}, \frac{q_2}{p_2}, \frac{q_3}{p_3})\) or \(M(-1,0; \frac{q'_1}{p'_1}, \frac{q'_2}{p'_2} )\), which contradicts to the definition of \(W\). If \(\partial M \neq \varnothing\) and \(C(h_1)\) is of type (e) or (f), then  by Remark \ref{paichu}, either \(M\) contains  singular fibers  with  slope-\(\frac{q}{2}\) where \(q\) is coprime to \(2\), or the base orbifold of \(M\)  is non-orientable. This also contradicts to the definition of \(W\). Therefore, \(C(h_1)\) must be  of type (d).

Suppose that \(n >1\),  then for each \(j=1, \cdots, n \), \(\partial (\widetilde{C}(h_{j})) \neq \varnothing\) and  each boundary component of  \(\widetilde{C}(h_{j})\) is incompressible in \(M\). By Lemma \ref{lem3},
 \(\widetilde{C}(h_{j})\) is a Seifert manifold.
It is a well-known fact that any incompressible torus in an irreducible Seifert-fibered manifold is isotopic to  
either a \emph{vertical} torus, i.e., a union of regular fibers, or \emph{horizontal} torus, i.e.,  transverse to all fibers
(see Hatcher \cite{Ha}). 
Since $M$ admits a unique Seifert fibering up to isotopy and $M$ is not homeomorphic to \( M(0,0;\frac{q_1}{p_1}, \frac{q_2}{p_2}, \frac{q_3}{p_3})\), 
 \(M\) does not contain any horizonal torus (by Corollary 3.12 in Jiang-Wang-Wu  \cite{JWW}  and Page 30 in Hatcher \cite{Ha}). Then we can choose a Seifert fibering of $M$ such that 
 each boundary component of  \(\widetilde{C}(h_{j})\) is a vertical torus of \(M\), which implies that \(\widetilde{C}(h_{j})\) admits a Seifert fibering induced by the Seifert fibering of \(M\).
Therefore, 
\(M\)
 is obtained from \(\widetilde{C}(h_{1}), \cdots, \widetilde{C}(h_{n}) \) due to the gluing  homeomorphisms that preserve the corresponding regular Seifert fibers.

By Lemma \ref{lem3}, \(\widetilde {C}(h_j)\) admits a unique Seifert fibering up to isotopy, unless \(\widetilde {C}(h_j)\) is homeomorphic to \(M (0,1;\frac{1}{2},\frac{1}{2})=M (-1,1;)\) or $T^2 \times I$ (see Corollary 3.12 in Jiang-Wang-Wu \cite{JWW}). 
If \(\widetilde{C}(h_{j})\) is homeomorphic to \(M (0,1;\frac{1}{2},\frac{1}{2})=M (-1,1;)\), then either \(M\) contains  singular fibers with slope-\(\frac{q}{2}\) where \(q\) is coprime to \(2\), or the base orbifold of \(M\)  is non-orientable. This contradicts the definition of \(W\). Therefore, 
 \(\widetilde {C}(h_{j})\) admits a unique Seifert fibering up to isotopy unless \(\widetilde {C}(h_j)\cong T^2 \times I\). 
Similar to the case that \(n=1\),  we can prove that  \(C(h_j)\) must be of type (d). The conclusions (1) and (2) of Proposition \ref{ITD SM} is proved.

Let  \(p_1, \cdots, p_s\) be the set of maximal non-separating points    of \(L\) (\(s=\beta_1 (L)\)), each of which corresponds to some torus in $\{\partial \widetilde {C}(h_{j})|j=1,\cdots, n \}$. Let $L'$ be the connected generalized graph obtained by cutting $L$ along these non-separating points, and $B$ be the base orbifold of $M$.
By conclusion (2), these non-separating points correspond to a set of non-separating circles \(c_1, \cdots, c_s\) of \(B\). 
Let \(\Sigma\) be the orbifold obtained  by cutting \(B\) along these circles, then  \(\Sigma\) is connected obviously. In fact, $\Sigma$ is obtained by gluing the base orbifolds of $  \widetilde {C}(h_{1}),\cdots, \widetilde {C}(h_{n}) $, where the gluing relationship is determined by $L'$.
Note that the genus of the base orbifold of each $\widetilde {C}(h_{j})$ is $0$. It is 
easy to observe that the genus of $\Sigma$ is $0$. Thus  the genus of \(B\) is equal to \(s=\beta_1(L)\). The proof of Proposition \ref{ITD SM} is completed.

\end{proo}

\section{Proof of Theorem \ref{t.realization}} \label{s.JSJ decompositions}

Let $W$ be an  ordinary graph manifold and  \(W=M_1 \cup \cdots \cup M_s\) be a  JSJ decomposition with the JSJ tori set \(\mathcal{T}\). 
In Section \ref{s.int}, we defined the indexed links related to the JSJ decompositions.
Let $l$ be an indexed link   related to the JSJ decomposition \(W=M_1 \cup \cdots \cup M_s\). Namely, $l$ satisfies the following condition:
\begin{enumerate}
\item \(l \cap (\cup_{T \in \mathcal{T}} T)=\varnothing\) and \(l\) contains both index-\(0\) knots and  index-\(2\) knots.
If \(T\in \mathcal{T}\) is  separating in \(W\), then  there is a knot of \(l\) with  index $0$ or  $2$
in 
each connected component of \( W| T\).

\item    For each \(i=1, \cdots ,s\), there is a Seifert fibering of $M_i$, such that
$l_i=l \cap M_i$ is a union of fibers which includes all of the singular fibers, and every singular fiber knot  is either index-$0$ or index-$2$. 
 
\item  
Let $x_i$ be the number of index-\(1\) knots of $l_i$, \(z_i\) be the number of other knots of $l_i$, \(b_i\) be the number of boundary components of \(M_i\),
 and \(g_i\) be the genus of the base orbifold of \(M_i\). Then \( z_i+b_i=x_i-2g_i+2 \).

\end{enumerate}

\begin{clai}
\label{c.rela}

Let \(m_i\) be the number of singular fibers of \(M_i\), then we have the following inequalities:
\begin{enumerate}
\item
 \(\Sigma_{i=1}^s z_i \geq 2\),
\(b_i +z_i \geq 2\), \(m_i \leq z_i\), and \(x_i \geq 1\).
\item
 If \(b_i=0\) or \(g_i \geq 1\), then \(x_i \geq 2\).
\item
If \(g_i=b_i=0\), then \(m_i \geq 4\).
\item
If \(g_i=0\) and \(b_i=1\), then \(m_i \geq 2\).
\end{enumerate}
\end{clai}
\begin{proo}
Since \(l\) contains both index-\(0\) knots and  index-\(2\) knots,
\(\Sigma_{i=1}^s z_i \geq 2\).  
If \(b_i =0\), then \(W \cong M_i\) and
 \(z_i \geq 2\). 
If \(b_i=1\), then \(\partial M_i\) corresponds to a separating JSJ torus \(T \in  \mathcal{T}\). By the condition (1) of the definition of \(l\),  there is a knot of \(l\) with  index $0$ or  $2$
in 
$M_i$, i.e.,
\(z_i \geq 1\). Thus, 
 we always have \(b_i +z_i \geq 2\). 
Since every singular fiber knot  is either index-$0$ or index-$2$,
 we have \(m_i \leq z_i\).

Suppose that \(g_i=0\), then \(b_i +z_i= x_i+2\).
If \(b_i =0\), then  \(M_i \cong W\). 
Since \(W\) is an ordinary graph manifold, \(W\) cannot be homeomorphic to  \(M (0,0;\frac{q_{1}}{p_{1}},\frac{q_{2}}{p_{2}}, \frac{q_{3}}{p_{3}})\). Then 
 \(z_i \geq m_i \geq 4\), which implies that  \(x_i \geq 2\). 
If \(b_i=1\) and \(m_i \leq 1\), then \(M_i \cong S^1 \times D^2\), which 
cannot be a Seifert piece of a closed graph manifold. Thus  when \(b_i=1\), we have \(z_i \geq m_i \geq 2\), which implies that \(x_i \geq 1\).
If \(b_i =2\) and \(z_i=0\), then \(M_i \cong T^2 \times I\). This  contradicts that  
$M_i$ admits a unique Seifert fibering up to isotopy.
Thus  when \(b_i=2\), we have \(z_i \geq 1\), then \(x_i \geq 1\). When \(b_i >2\), by the 
 equation 
\(b_i +z_i= x_i+2\), we have \(x_i \geq 1\).

When \(g_i \geq 1\), we have \(2 \leq b_i +z_i= x_i-2g_i+2 \leq x_i\).
Therefore, we always have \(x_i \geq 1\). In particular, if \(b_i=0\) or \(g_i \geq 1\), then \(x_i \geq 2\).

\end{proo}

Let \(v_1, \cdots, v_s\) be \(s\) vertices in \(\mathbb{R}^3\). If \(T\in \mathcal{T}\) is adjacent  to \(M_i\) and \(M_j\), then we connect an edge ending at $v_i$ and $v_j$, where \(i,j \in \{1, \cdots, s\}\) and it is possible that \(i=j\).
According to the JSJ decomposition \(W=M_1 \cup \cdots \cup M_s\), we get a $1$-complex \(G\),  called a \emph{JSJ graph} of $W$.

 Let \( \mbox{Star}(v_i)\) be a small neighborhood of \(v_i\) in \(G\), then $\mbox{Star}(v_i)$ has $b_i$ edges.   If
\( \mathcal{T}\) consists of \(r \geq 1\) tori, then
 \(G\) has \(r\) edges, and we denote  the edges of \(G\) by \(e_1, \cdots, e_r\). We denote by \(\mathcal{E}\) the set of the edges of \(G\) corresponding to separating JSJ tori.
We label the edge in  \( \mbox{Star}(v_i)\) by \(e_j\) if this edge lies in \(e_j\).

\begin{lemm}\label{c.conLya}
There is a connected graph \(L\) and
a surjective projection $\pi: L\to G$ such that:
\begin{enumerate}
\item
For each edge $e_j$ and each vertex \(v_i\) of \(G\), 
 $\pi^{-1} (e_j)$ is an edge of $L$, and  $L_i= \pi^{-1} (\mbox{Star}(v_i))$ is a connected generalized graph with \(b_i\) ends.

\item
The vertices of  $L_i$ consist of  $z_i$ degree $1$ vertices and \(x_i\) degree $3$ vertices,   and $\beta_{1}(L_i)=g_i$. Let  $y_i$ be the number of the vertices of $L_i$ adjacent to two degree $1$ vertices,
then $2y_i \leq m_i$.

\item  
If  
 \(v_{i}\) is  in a circle of \(G\),  then the midpoint of the edge in \(L_{i}\) 
  is non-separating in \(L\) unless this edge connects degree $1$ vertices or intersects \(\cup_{e \in \mathcal{E} } \pi^{-1}( e)\).

\end{enumerate}

\end{lemm}
\begin{proo}
For each  \( i=1, \cdots ,s\), we suppose that $L_i$ is a connected generalized graph with \(b_i\) ends.
Note that 
\( \mbox{Star}(v_i)\) has \(b_i\) edges. Then
 there is a one-to-one correspondence between the set of the edges in \( \mbox{Star}(v_i)\) and the set of the ends of \(L_i\). We label the ends of \(L_i\)  by the labels on the corresponding edges in \( \mbox{Star}(v_i)\). 

Suppose that \(g_i=0\), then \( b_i +z_i= x_i+2\).
If \(b_i=0\), then by Claim \ref{c.rela}, we have
 \(z_i \geq m_i \geq 4\), and
 \(L_i\) can be constructed as Figure \ref{Fgsec6.1} (a). 
If \(b_i=1\), then by Claim \ref{c.rela}, we have \(z_i \geq m_i \geq 2\), and \(L_i\) can be constructed as Figure  \ref{Fgsec6.1} (b).
If \(b_i \geq2\), then
 \(L_i\) can be constructed as  Figure  \ref{Fgsec6.1} (c). 
 
\begin{figure}[htbp]
\centering
\subfigure[ ]{\includegraphics[width=0.2\textwidth]{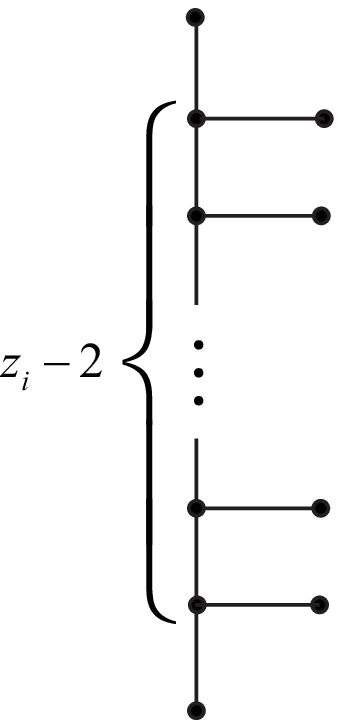}}
\hspace{.60in}
\subfigure[ ]{\includegraphics[width=0.2\textwidth]{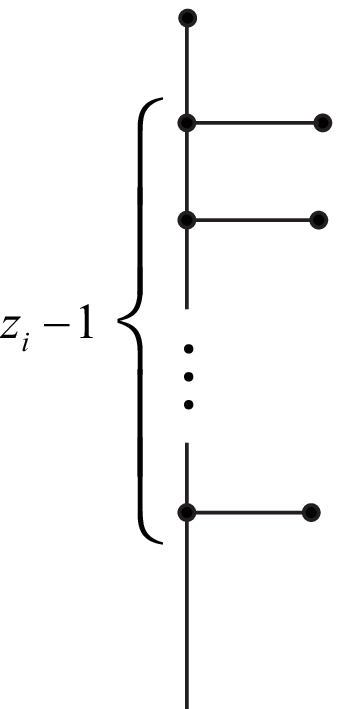}}
\hspace{.60in}
\subfigure[ ]{\includegraphics[width=0.2\textwidth]{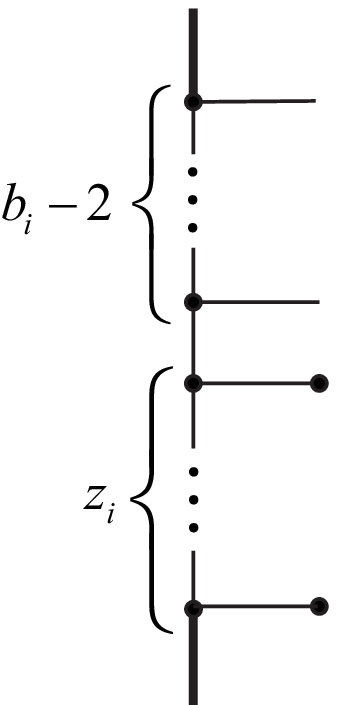}}

\caption{These three generalized graphs are trees. There are two vertices that are adjacent to two degree $1$ vertices in Figure (a), there is one vertex adjacent to two degree $1$ vertices in Figure (b), and there is no vertex adjacent to two degree $1$ vertices in Figure (c). }
\label{Fgsec6.1}
\end{figure}

Suppose that \(g_i \geq 1\). Since \(2 \leq b_i +z_i= x_i-2g_i+2\), we can  construct \(L_i\) as shown in Figure  \ref{Fgsec6.2}. It is easy to observe that the above \(L_i\) that we constructed  satisfy the condition (2). 
\begin{figure}[htbp]
\centering
\includegraphics[scale=0.6]{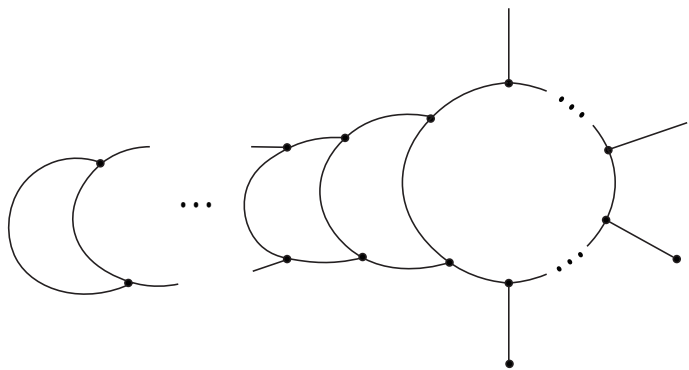}
\caption{This generalized graph  supports \(z_i\) degree $1$ vertices and \(b_i\) ends, and its first Betti number is \(g_i\).}
\label{Fgsec6.2}
\end{figure}

By
gluing \(L_1 , \cdots, L_s\) along the ends labeled by the same letters,
we can get a connected graph \(L\). 
Let \(L'_i\) be the connected generalized graph obtained from \(L_i\) by removing the ends and the edges connecting ends. By Claim \ref{c.rela}, \(x_i \geq 1\), then \(L'_i \neq \varnothing\).
By mapping each $L'_i$ into a vertex, we can get \(G\) from \(L\). Thus
we can construct 
  a surjective projection \(\pi: L \to G\) such that 
 \( \pi^{-1} (\mbox{Star}(v_i))\) is equal to \(L_i\) for each \(i=1, \cdots, s\), and   $\pi^{-1} (e_j)$ is an edge  of \(L\) for each edge $e_j$ of \(G\).
 Moreover, for any \(e \in \mathcal{E}\), the midpoint of \(\pi^{-1}( e)\) is separating in \(L\).

Suppose that \(v_i\) is in a circle of \(G\). 
If \(g_i \geq 1\), then the midpoint of the edge in \(L_{i}\) 
  is non-separating in \(L\) unless this edge connects degree $1$ vertices or intersects \(\cup_{e \in \mathcal{E} } \pi^{-1}( e)\) (see Figure \ref{Fgsec6.2}). If \(g_i=0\), then we choose a circle \(C_i\) of \(G\)   that contains \(v_i \). Obviously,
\( \mbox{Star}(v_i) \cap C_i\) contains two edges, then we make them  correspond to the 
  two bold edges  in  Figure  \ref{Fgsec6.1} (c). Thus the conclusion (3) in this case is satisfied.
 Based on the above construction, the surjective projection \(\pi\) satisfies the conclusions of Lemma \ref{c.conLya}.

\end{proo}

Let 
\(K_{i}\) be the set of the index-\(0\) or index-\(2\) knots of \(l_i\), and  \(K^{0}_{i} \subset K_{i}\) be the set of the  singular fibers of \(M_i\). Let \(S_i\) be the set of the degree $1$ vertices of \(L_i\),
and \(S'_i\) be the set of the vertices of \(L_i\)  adjacent to two degree $1$ vertices.
Since \(L_i\) contains \(z_i\) degree $1$ vertices and  \(2y_i \leq m_i \leq z_i\), we can construct a  bijection \(\sigma_i : S_i \to K_i\) such that each degree $1$ vertex adjacent to a vertex in \(S'_i\) is contained in \(\sigma^{-1}_i (K^{0}_i)\).

We endow the orientation on the edges connecting a degree $1$ vertex in \(L_i\) such that the degree $1$ vertices mapping into index-\(0\) knots under \(\sigma_{i}\) are sink vertices, and  degree $1$ vertices mapping into index-\(2\) knots under \(\sigma_{i}\) are source vertices.
Since \(l\) contains both index-\(0\) knots and  index-\(2\) knots, \(L\) must contain both source degree $1$ vertices and sink degree $1$ vertices.

\begin{lemm}
 \(L \in \mathcal{S}\), where \(\mathcal{S}\) is defined in Definition \ref{d.graphori}.

\end{lemm}
\begin{proo}
In fact,
we only need to prove that for any separating point  of \(L\), each connected component obtained by cutting \(L\) along this point
 must contain degree $1$ vertices.
Now we prove it by contradiction.
Suppose that there is a separating point \(p\) of \(L\) such that a connected component \(L'\) of \( L\smallsetminus \{p\}\) does not contain any degree $1$ vertex. 
Obviously, \(L'\) must contain circles, then $\beta_{1}(L')\geq 1$.  Since \(L\) contains at least two degree $1$ vertices, the edge \(e'\) containing \(p\)   cannot  connect degree $1$ vertices. 

If  \(\pi(L'\cup e') \cap \mathcal{E} \neq \varnothing\), then 
 there is a point \(p_0\) in \(L'\cup e'\) such that \(\pi(p_0)\) is in an edge in  \(\mathcal{E}\). According to the construction of \(\pi\),
 \(p_0\) is separating  in \(L\) and there is a component of \(L\smallsetminus \{p_0\}\)  contained in \(L'\cup e'\).
 By the condition (1) of  the definition of the indexed link \(l\), each component of \(L\smallsetminus \{p_0\}\) must contain degree $1$ vertices, which implies that 
\(L'\) must contain degree $1$ vertices. This contradicts to the supposition,
thus
 \(\pi(L'\cup e') \cap \mathcal{E} = \varnothing\). 

If \(\mathcal{E}= \varnothing\),
then either  each vertex \(v_i\) of \(G\) is in a circle of \(G\),  or \(G\) is a vertex.
For the first case, by the conclusion (3) of Lemma \ref{c.conLya}, it is easy to observe that the midpoint of the edge in \(L \)  is non-separating in \(L\) unless this edge connects degree $1$ vertices. 
For the second case, we may assume that \(G\) is the vertex \(v_1\), then \(L=L_1\). Since $\beta_{1}(L')\geq 1$,  we have \(\beta_{1}(L) \geq 1\). 
According to the construction of \(L\) in Lemma \ref{c.conLya}, the midpoint of the edge in \(L\)  is non-separating in \(L\) unless this edge connects degree $1$ vertices.  However, \(e'\) cannot  connect degree $1$ vertices, which  contradicts that \(p\) is  separating  in \(L\).

If \(\mathcal{E} \neq \varnothing\), then  by cutting \(G\) along the midpoint of each edge in \(\mathcal{E}\), we can get
a connected component \(G_0\)  containing \(\pi (L'\cup e' ) \).
This is because that \(\pi(L'\cup e') \cap \mathcal{E} = \varnothing\). Then either each vertex   of \(G_0\) is in a circle of \(G_0\),  or \(G_0\) consists of  a vertex and \(G_0\cap \mathcal{E}\).
Let \(L_0= \pi^{-1}(G_0)\), then 
 \((L'\cup e' ) \subset L_0\). 
For the first case, by the conclusion (3) of Lemma \ref{c.conLya},  it is easy to observe that the midpoint of the edge in \(L_0 \)  is non-separating in \(L\) unless this edge connects degree $1$ vertices or intersects  \(\cup_{e \in \mathcal{E} } \pi^{-1}( e)\). 
For the second case,  we may assume that \(G_0\) contains the vertex \(v_1\), then
\(G_0\) is a neighborhood of \(v_1\) in \(G\), which implies that 
 \(L_0=L_1\) (regardless of the length of the edges). Since  \(\beta_{1}(L') \geq 1\), we have \(g_1=\beta_{1}(L_1) \geq 1\).
According to the construction of \(L_1\) in Lemma \ref{c.conLya}, the midpoint of the edge in \(L_1\)  is non-separating in \(L\) unless this edge connects degree $1$ vertices or intersects  \(\cup_{e \in \mathcal{E} } \pi^{-1}( e)\). However, \(e'\) cannot  connect degree $1$ vertices and \(e' \cap (\cup_{e \in \mathcal{E} } \pi^{-1}( e)) = \varnothing\), 
 which contradicts that \(p\) is   separating   in \(L\). Therefore, \(L \in \mathcal{S}\).

\end{proo}

By Lemma \ref{orientation}, we can endow  orientations on the other edges of  \(L\),
so that \(L\) is an abstract Lyapunov graph,  each of the source and  sink vertices is a degree $1$ vertex, and each saddle vertex is a degree $3$ vertex. 

\textbf{Proof of Theorem \ref{t.realization}:}
%
%
%
%
%
Let  \( i \in \{1, \cdots ,s\}\).
It is note that  \(L_i\)  has been oriented to  an abstract Lyapunov graph containing \(x_i\) saddle vertices. 

Let \(  \mathcal{X}_i\) be the set of the saddle vertices of \(L_i\).
 By Claim \ref{c.rela},
 \(x_i \geq 1\), and when \(b_i=0\), \(x_i \geq 2\). Thus
 \(  \mathcal{X}_i \neq \varnothing\) and
each saddle vertex is adjacent to at most two degree $1$ vertices.
  From now on, we give the corresponding relationship between the block \(L_v\)  associated to a saddle vertex \(v\) of \(L_i\) and a \(\widetilde{C}(h_v)\) in Lemma \ref{lem3}, where \(C(h_v)\) is a  \(1\)-FRH of type (d) in Lemma \ref{lem2}. 
The block \(L_v\)  associated to  \(v\)  is defined in Section \ref{ss.Lya}.

If \(v\) is adjacent to  two degree $1$ vertices \(v'\) and \(v''\), then \(\sigma_i (v'),  \sigma_i (v'')\in K_i^0\). Namely,
 \(\sigma_i (v')\) and \(\sigma_i (v'')\) are two singular fibers of \(M_i\). Suppose that \(\frac{q'}{p'}\) and \(\frac{q''}{p''}\) are the slopes of  \(\sigma_i (v')\) and \(\sigma_i (v'')\) respectively. Then
we make \(L_v\) correspond to \(\widetilde{C}(h_v) \cong M (0,1;\frac{q'}{p'},\frac{q''}{p''})\).
 If \(v\) is adjacent to  one degree $1$ vertex \(v'\) and \(\sigma_i (v') \in K_i^0\), then 
we make \(L_v\) correspond to \(\widetilde{C}(h_v)\cong M (0,2;\frac{q'}{p'})\), where
\(\frac{q'}{p'}\) is the slope of  \(\sigma_i (v')\).
If \(v\) is adjacent to  one degree $1$ vertex \(v'\) and \(\sigma_i (v') \notin K_i^0\), then 
we make \(L_v\) correspond to \(\widetilde{C}(h_v)\cong M (0,2;)\).
 In addition, if \(v\) is not adjacent to  any degree $1$ vertex, we make \(L_v\) correspond to \(\widetilde{C}(h_v) \cong M (0,3;)\). 

Note that the above Seifert fibering of \(\widetilde{C}(h_v)\) admits the cores of the 
RHs associated to natural FRH decomposition of \(\widetilde{C}(h_v)\) to act  as fibers.
Moreover, the cores of the RHs consists of all singular fibers and some regular fibers of 
\(\widetilde{C}(h_v)\).
By Theorem \ref{t.thm1}, we can construct an NMS flow on \(\widetilde {C}(h_v)\)  such that the natural FRH decomposition of \(\widetilde {C}(h_v)\) is a FRH decomposition of this flow.
In fact, we can choose the above  \(\widetilde {C}(h_v)\) good enough  such that
\(L_v\) is a Lyapunov graph of this flow and 
  corresponds to this natural FRH decomposition. 


By Lemma \ref{lemorder}, we get a gluing order of the blocks associated to  saddle vertices of \(L_i\), then we can   get the corresponding gluing order of  \(\{\widetilde {C}(h_v) |  v \in \mathcal{X}_i\}\).
Since
 \(\sigma_i : S_i \to K_i\) is a bijection and \(K_i^0 \subset K_i\), the  slopes of the singular fibers in \(\sqcup _{  v \in \mathcal{X}_i} \widetilde {C}(h_v)  \) are the slopes of the singular fibers of \(M_i\) exactly. Therefore,
by gluing \(\{\widetilde {C}(h_v) |  v \in \mathcal{X}_i\}\) with a series of homeomorphisms preserving fibers, we can get a Seifert manifold \(N_i= M (g'_i,b_i;\frac{q_{1}}{p_{1}},\frac{q_{2}}{p_{2}},\cdots ,\frac{q_{m_i}}{p_{m_i}})\),  where \(\frac{q_{1}}{p_{1}},\frac{q_{2}}{p_{2}},\cdots ,\frac{q_{m_i}}{p_{m_i}}\) are the slopes of the singular fibers of \(M_i\) exactly.

Based on the above discussion, we get a FRH decomposition \(N_i= (\partial_{-} N_i \times I) \cup _{  v \in \mathcal{X}_i} \widetilde {C}(h_v)\) of \(N_i\). By Theorem \ref{t.thm1}, we can construct an NMS flow \(\psi^i_t\) on \(N_i\), such that the decomposition \(N_i= (\partial_{-} N_i \times I)\cup _{  v \in \mathcal{X}_i} \widetilde {C}(h_v)\) is a FRH decomposition of \(\psi_t^i\), and
$L_i$  is the Lyapunov graph of  \(\psi^i_t\)  
 corresponding to this decomposition. 
 Similar to the proof of Proposition \ref{ITD SM}, we have \(g'_i = \beta_{1}(L_i)=g_i\), which implies that \(N_i \cong M_i\).

 Due to the surjective projection $\pi: L\to G$,  we can get a FRH decomposition \(W = \cup_k C_k\) of an NMS flow \(\phi_t\) on  \(W\) by suitably gluing the above FRH  decompositions of \(\psi^1_t, \cdots, \psi^s_t\), so that $L$  is the Lyapunov graph of $\phi_t$ 
corresponding to the  decomposition \(W = \cup_k C_k\). 
Note that the intersection of indexed link \(l'\) of \(\phi_t\) and each RH  associated to the decomposition \(W=  \cup_{k} C_{k}\) is a core of this RH.
It is easy to observe that \(l'\)  is isotopic to \(l\)  and the index of each knot in \(l'\) is the same as the index of the corresponding  knot in \(l\).
Thus we can perturb \(\phi_t\) such that \(l\) is its indexed link.
Theorem \ref{t.realization} is proved.
 
%
%

\begin{rema}\label{r.easyconstruction}
\begin{enumerate}
\item According to the construction of $\phi_t$, it is easy to prove that
\(\phi_t\) is related to  incompressible torus decompositions.
\item
In this section, we construct a FRH decomposition of some NMS flow on $W$ such that the indexed link of this flow is   $l$. In fact,
for this purpose, 
 we only need to construct a connected graph $L$ that satisfies the conditions (1), (2) of Lemma \ref{c.conLya} and the following conditions.
\begin{itemize}
\item Each connected component obtained by cutting $L$ along any separating point must contain degree $1$ vertices.
\end{itemize}
 \end{enumerate}
\end{rema}

\section{Proof of Theorem \ref{t.main0}}\label{s.change regular fibers}

Let $W$ be an  ordinary graph manifold and  \(W=M_1 \cup \cdots \cup M_s\) be a  JSJ decomposition with the JSJ tori set \(\mathcal{T}\).  
Let $l$ be an indexed link   related to the JSJ decomposition \(W=M_1 \cup \cdots \cup M_s\) (defined in Section \ref{s.int}).

An incompressible torus set $\mathfrak{T}$  of $W$ \emph{is  related to} $l$ if:
\begin{enumerate}
\item \(\mathcal{T}\subset \mathfrak{T}\), and $l\cap (\cup_{T \in \mathfrak{T}} T)=\varnothing$.
\item If $T \in  \mathfrak{T}$ is separating in $W$, then 
there is a knot of $l$ with index $0$ or $2$ in 
each component of $W|T$.
\item   each connected component $C$ of $W | \mathfrak{T}$ is atoroidal.

\item $l\cap C$  contains only one  index-$1$ knot  and the number of knots in $l\cap C$ is equel to $4-n$, where $n$ is the number of  the components of $\partial C$.


\end{enumerate}

\begin{rema}
\begin{enumerate}
\item
As the discussion in Section \ref{s.JSJ decompositions}, there exists an incompressible torus set $\mathfrak{T}$  related to $l$.
\item By Remark \ref{r.easyconstruction}, the incompressible torus set related to $l$ is easy to construct.

\end{enumerate}
\end{rema}

\begin{figure}[htbp]
\centering
\subfigure[ ]{\includegraphics[width=0.55\textwidth]{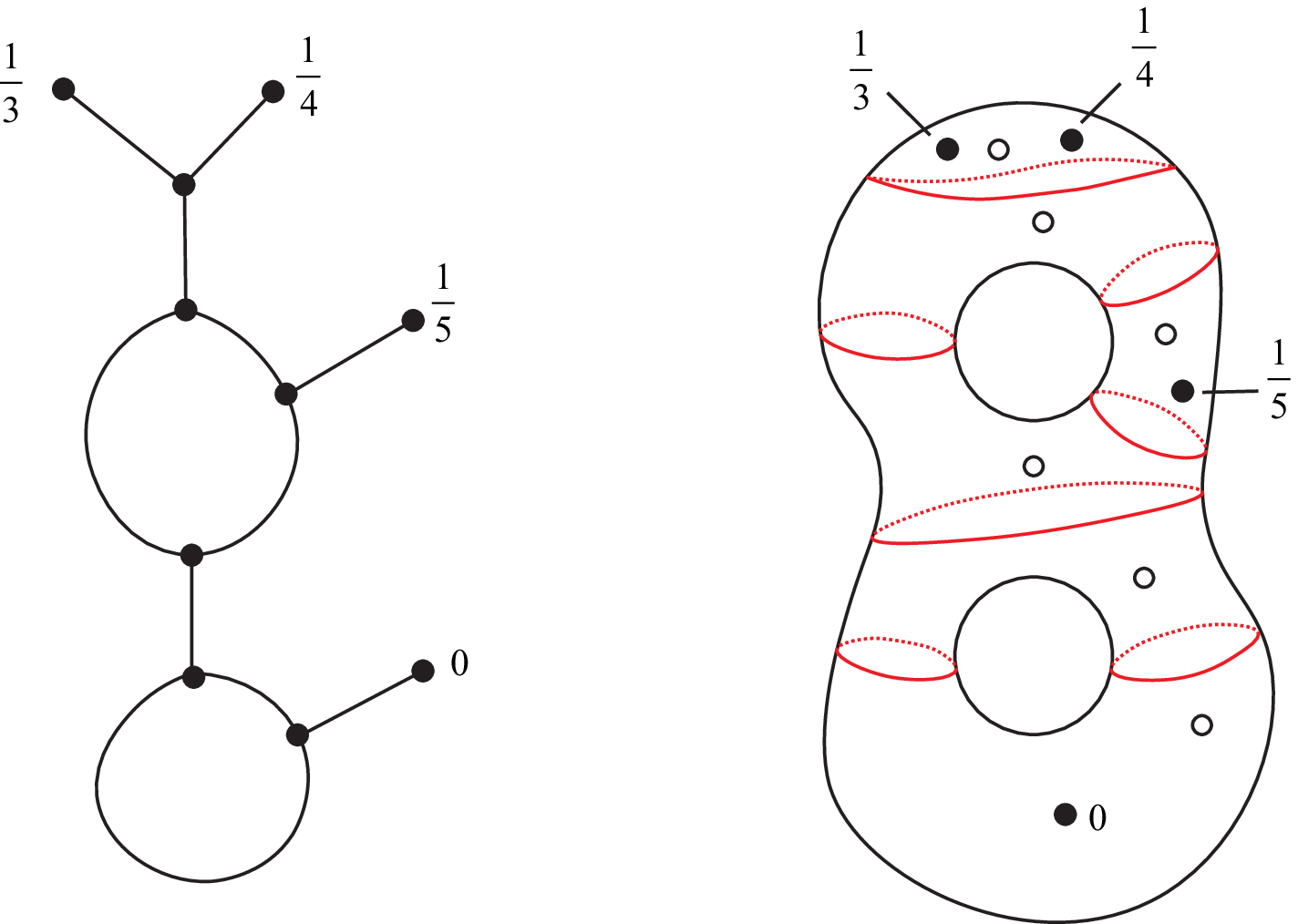}}
\hspace{.60in}
\subfigure[ ]{\includegraphics[width=0.55\textwidth]{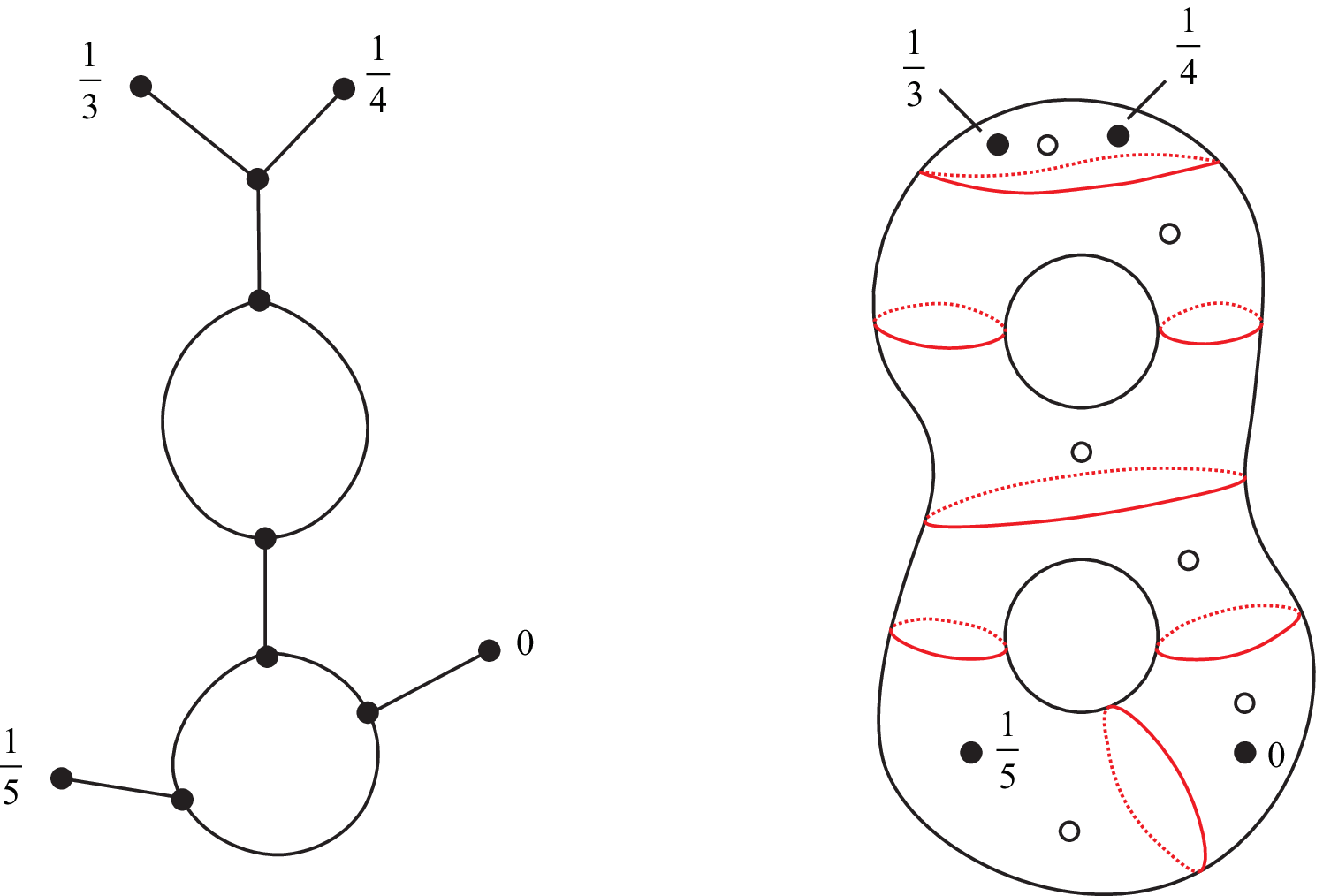}}
\caption{The number is the slope of the corresponding fiber.
On the right side of each subgraph is the base orbifold of $M$, where   hollow circles correspond  to index-$1$ knots and  solid circles correspond  to other knots of $l$. }
\label{Fgpartation}
\end{figure}

%
%
%
%
\begin{example}
Let $M=M(2,0; \frac{1}{3}, \frac{1}{4}, \frac{1}{5})$ be a Seifert manifold, and $l$ be an indexed link related to the JSJ decomposition of $M$. Suppose that $l$ has $6$ index-$1$ knots. Then the circles in  Figure \ref{Fgpartation} (a) and the circles in  Figure \ref{Fgpartation} (b) are corresponding to two incompressible torus sets related to $l$, respectively.

\end{example}

\textbf{Operation of changing regular fibers: }
Choose an  incompressible torus set $\mathfrak{T}$  related to $l$. For  each component $C$ of $W | \mathfrak{T}$  homeomorphic to $T^2 \times I$, $l\cap C=k_1\cup  k_2$ are two parallel regular fibers of  a Seifert fibering of $C$. 
We replace $ k_1,  k_2$ with two parallel  torus knots $c_1, c_2$ in $C$, where $\text{Ind}(k_i)=\text{Ind}(c_i)$ for $i=1,2$.

\begin{rema}
 ``Operation of changing regular fibers" depends on ``incompressible torus set related to $l$", i.e., two different incompressible torus sets related to $l$  maybe induce  two different indexed links by applying Operation  of changing regular fibers to $l$.

\end{rema}

\begin{theorem}\label{t.incindlink1}
Let \(l\) be an indexed link in an ordinary graph manifold $W$. Then \(l\) is the indexed link 
of an NMS flow on $W$ related to  incompressible torus decompositions if and only if
there is an indexed link $l'$ related to 
a JSJ decomposition of $W$ such that $l$
can be obtained from  $l'$ by  applying at most one step of  Operation of changing regular fibers.
\end{theorem}

Notice that applying at most one step of  Operation of changing regular fibers 
 actually means changing a finite number of regular fiber pairs.

\begin{proo}
\emph{Necessity}. 
Let   \(\phi_t\) be an NMS flow on $W$ related to   incompressible torus decompositions, and 
$l$ be the indexed link of  \(\phi_t\). 

Since \(W\) is a closed manifold,
any NMS flow on \(W\) must contain both attracting closed orbits and  repelling closed orbits. Namely, \(l\) must contain both  index-\(0\) knots and  index-\(2\) knots.
By Proposition \ref{prop2}, there     is a JSJ decomposition \(W= M_1 \cup \cdots \cup M_s\) 
 with the JSJ tori set \(\mathcal{T}\) such that
\(\phi_t\) is transverse to each 
  \(T \in \mathcal{T}\).
Then \(l \cap T= \varnothing\).
If 
\(T'\in \mathcal{T}\) is separating in \(W\), then we suppose that \(W| T'= W' \sqcup W''\). 
Since \(\phi_t\) is transverse to \(T'\) and \(\partial W'\) is connected, 
\(\phi_t| _{W'}\) 
must contain  attracting closed orbits or  repelling closed orbits. 
 Then  \(l \cap W'\)
must contain at least one knot with index $0$ or $2$. Similarly, \(l \cap W''\) must also contain at least one knot with index $0$ or $2$.

For each \(i=1, \cdots, s\),
let 
\(l_i =l \cap M_i\), and we suppose that  \(l_i\) consists of  \(x_i\) index-\(1\) knots and \(z_i\) other knots.
 Proposition \ref{prop2} shows that 
\(\phi_t |_{M_i}\) is an NMS flow related to   incompressible torus decompositions, then \(x_i \geq1\). Let  \(M_i= (\partial _{-}M_i \times I) \cup_{j=1}^{x_i}\widetilde{C}(h_j ^i)\) be an incompressible torus decomposition  of \(\phi_t|_{M_i}\).

By Proposition \ref{ITD SM}, 
each \(C(h_{j}^i)\) is of type (d) in Lemma \ref{lem2}, and  $\widetilde {C}(h_j^i)$ is a Seifert manifold with incompressible boundary.
Note that the intersection of \(l_i\) and each RH   associated to \(\widetilde{C}(h_{j}^i)\) is a core of this RH.
According to the proof of Lemma \ref{lem3}, \(l_i \cap \widetilde{C}(h_j^i)\) consists of all the singular fibers and some regular fibers of a Seifert fibering \(\widetilde{C}(h_j^i)\),  where  each singular fiber knot  is either index-$0$ or index-$2$. 
For the convenience of description, we refer to the above Seifert fibering of  \(\widetilde{C}(h_j^i)\)  as the \emph{natural Seifert fibering} of \(\widetilde{C}(h_j^i)\).
Proposition \ref{ITD SM} illustrates that
 \(M_i\) is obtained from \(\left \{ \widetilde {C}(h_j^i)|j=1,\cdots, x_i \right \}\)  by the gluing  homeomorphisms that preserve the corresponding regular Seifert fibers.  
Then we   choose  a Seifert fibering of $M_i$, such that $\partial \widetilde {C}(h_j^i)$ is a   union of vertical tori for each $j$.

\textbf{Case 1}. For any $i\in \{1, \cdots, s\}$ and any $j\in\{1, \cdots, x_i\}$, the natural Seifert fibering of \(\widetilde{C}(h_j^i)\) is isotopic to the restriction of the Seifert fibering of $M_i$ to \(\widetilde{C}(h_j^i)\).

Then there is a Seifert fibering of $M_i$, such that 
\(l_i\) is a union of fibers which includes all of the singular fibers, and each singular fiber knot  is either index-$0$ or index-$2$. 
 Let \(b_i\) be the number of boundary components of \(M_i\), and \(g_i\) be the genus of the base orbifold of \(M_i\).
Let  \(L_i\) be the  Lyapunov graph of  \(\phi_t |_{M_i}\), which  corresponds to the  decomposition \(M_i= (\partial _{-}M_i \times I) \cup_{j=1}^{x_i}\widetilde{C}(h_j^i)\).
 By Proposition \ref{p.ends}, each of the source and sink vertices is a degree $1$ vertex. Then the vertices of \(L_i\) consists of \(x_i\) saddle vertices and \(z_i\) degree $1$ vertices, and $L_i$ has \(b_i\) ends.
 By Proposition \ref{ITD SM}, \(\beta_1 (L_i)=g_i\), and  each saddle  vertex of \(L_i\) connects three edges.

Let  \(p^{i}_1, \cdots, p^{i}_{g_i}\) be a set of maximal non-separating points    of \(L_i\). By cutting \(L_i\) along these points, we can get a tree \(L'_i\) that has \(x_i\) degree $3$ vertices, \(z_i\) degree $1$ vertices and \(b_i+2g_i\) ends. By induction on the number of degree $3$ vertices in \(L'_i\), it is easy to prove that \(z_i+b_i+2g_i=x_i+2\). Namely,  \(z_i+b_i=x_i-2g_i+2\). Therefore $l$ is an indexed link related to  the JSJ decomposition  \(W= M_1 \cup \cdots \cup M_s\).

\textbf{Case 2}. There exists   $i'\in \{1, \cdots, s\}$ and   $j'\in\{1, \cdots, x_i\}$ such that the natural Seifert fibering of \(\widetilde{C}(h_{j'}^{i'})\) is not isotopic to the restriction of the Seifert fibering of $M_i$ to \(\widetilde{C}(h_{j'}^{i'})\).

By Proposition \ref{ITD SM},   \(\widetilde{C}(h_{j'}^{i'}) \cong T^2 \times I\).
By Lemma \ref{lem3}, $l_{i'}\cap \widetilde{C}(h_{j'}^{i'})$ consists of an index-$1$ knot $c_1$ and a knot $c_2$ that is not of index-$1$. Choose two regular fibers $k_1, k_2$ in the restriction of the Seifert fibering of $M_i$ to \(\widetilde{C}(h_{j'}^{i'})\), and endow the indices on $k_1, k_2$ such that $\text{Ind}(k_1)=\text{Ind}(c_1)$ and $\text{Ind}(k_2)=\text{Ind}(c_2)$.
Let $R$ be a solid torus in $\widetilde{C}(h_{j'}^{i'})$ such that $k_2$ is a core of $R$, $\partial R$ is  vertical associated to the restriction of the Seifert fibering of $M_i$ to \(\widetilde{C}(h_{j'}^{i'})\),
 and $R \cap k_1 = \varnothing$. Then $\overline{\widetilde{C}(h_{j'}^{i'})\smallsetminus R}$ is homeomorphic to the circle bundle over a pair-of-pants. Therefore, we get a new FRH decomposition of $\widetilde{C}(h_{j'}^{i'})$, such that the cores of the RHs  
are $k_1, k_2$. Similarly, we discuss all the  \(\widetilde{C}(h_j^i)\) that are homeomorphic to $T^2 \times I$. Finally, we get a  FRH decomposition of a new NMS flow $\phi_t'$ on  $W$. 
Let $l'$ be the indexed link of $\phi'_t$. Similar to the discussion in Case 1, $l'$  is  an indexed link related to the JSJ decomposition  \(W= M_1 \cup \cdots \cup M_s\).
Moreover, $l$ is obtained from $l'$ by applying Operation of changing regular fibers.  
%
\\
\par \emph{Sufficiency}.  Let \(W=M_1 \cup \cdots \cup M_s\) be a  JSJ decomposition of $W$, and
 $l'$ be an indexed link   related to the JSJ decomposition \(W=M_1 \cup \cdots \cup M_s\). According to the discussion in Section \ref{s.JSJ decompositions}, $l'$ is the indexed link of an NMS flow on $W$ related to incompressible torus decompositions.

Let  $\mathfrak{T}$ be an incompressible torus set  related to $l'$. Suppose that $l$ is the indexed link obtained from $l'$ by Operation of changing regural fibers along $\mathfrak{T}$. From now on, we prove that \(l\) is the indexed link 
of some NMS flow on $W$ related to  incompressible torus decompositions.

Since $M$ admits a unique Seifert fibering up to isotopy and $M$ is not homeomorphic to \( M(0,0;\frac{q_1}{p_1}, \frac{q_2}{p_2}, \frac{q_3}{p_3})\), 
 \(M_i\) does not contain any horizonal torus for each $i=1, \cdots, s$. Then there is a Seifert fibering of $M_i$ such that each torus $T \in \mathfrak{T}\cap M_i$  is a vertical torus and $l'\cap M_i$ are fibers. 
Let $C$ be a  connected component  of $W | \mathfrak{T}$.
Then $C$ is a Seifert manifold with the induced Seifert fibering. 
By the definition of $l'$, $l' \cap C$ consists of  all the singular fibers and some regular fibers of \(C\), where each singular fiber knot  is either index-$0$ or index-$2$. 

Recall that  $C$ is atoroidal and contains only one  index-$1$ knot of $l$. Moreover, 
the number of knots in $l'\cap C$ is equel to $4-n$, where $n$ is the number of  the components of $\partial C$.
Then there are three possibilities.
\begin{itemize}
\item  $C\cong M (0,1;\frac{q_1}{p_1},\frac{q_2}{p_2})$ ($p_1 \cdot p_2 \neq 1$). Then   $l'\cap C$ consists of 
two singular fibers and one regular fiber, where  the index of regular fiber is $1$.

\item $C\cong M (0,2;\frac{q}{p})$. If $p\neq 1$, then  $l'\cap C$ consists of the singular fiber and one regular fiber  with index $1$. Otherwise, $l'\cap C$ consists of two regular fibers, where one knot is  index-$1$ and the other knot is  either index-$0$ or index-$2$.

\item $C\cong M (0,3;)$. Then   $l'\cap C$ is a regular fiber   with index $1$.
\end{itemize}

Let $N$ be the number of the knots in $l'$, and  $V$ be the union of 
the different vertices $v_1, \cdots, v_N$ in \(\mathbb{R}^3\). 
 Choose a one-to-one correspondence between the  vertices in $V$ and the knots in $l'$. We label the vertex by the index of the corresponding knot. 
Let $r$ be the number of the vertices labeled by $1$.
Assume that 
$v_1, \cdots, v_r$ are labeled by $1$. 
For \(j,k \in \{1, \cdots, r\}\), we connect an edge ending at $v_j$ and $v_k$ if  their corresponding knots are in adjacent components of $W|\mathfrak{T}$. For each vertex $v_j$ labeled by $1$, we connect an edge between $v_j$ and the vertex whose  corresponding knot is in the same components of $W|\mathfrak{T}$ as the knot corresponding $v_j$. 
Then we get a connected graph $L$ such that each vertex either connects one edge or three edges.
We endow the orientation on the edges connecting  the vertices labeled by \(0\) or $2$ in \(L\) such that the vertices labeled by \(0\)   are sink vertices  and  the vertices labeled by \(2\)  are source  vertices.  
 It is easy to observe that $L \in \mathcal S$ (defined in Section \ref{ss.Lya})

By Lemma \ref{orientation}, we can endow  orientations on other edges of  \(L\),
such that \(L\) is an abstract Lyapunov graph where  the vertex connecting three edges is a saddle vertex.  For each  component $C$ of $W|\mathfrak{T}$, let $L_c$ be the block   associated to  the saddle vertex that corresponds the index-$1$ knot of $l'\cap C$. 
Using Lemma \ref{lem3}, we can construct a FRH decomposition of some NMS flow on  $C$  such that 
\begin{itemize}
\item
the $1$-FRH in this decomposition is of type (d) in Lemma \ref{lem2};
\item
$l'\cap C$ consists of the cores of the RHs associated to this decomposition;
 \item \(L_c\) is a Lyapunov graph of this NMS flow and 
 corresponds to this  FRH decomposition. 
\end{itemize}
 
By Lemma \ref{lemorder} and Theorem \ref{t.thm1}, there is a FRH decomposition of some NMS flow $\phi'_t$ on $W$ such that $L$ is a Lyapunov graph of $\phi'_t$ and $l'$ is the indexed link of $\phi'_t$.
For each component $C'\cong T^2 \times I$ of $W|\mathfrak{T}$, we can adjust suitably the FRH decompositon of $C'$, and then we  construct a   FRH decomposition of a new NMS flow $\phi_t$ on $W$, such that $l$ is the indexed link of $\phi_t$. By Lemma \ref{lemITD}, $\phi_t$ is related  to incompressible torus decompositions.
Theorem \ref{t.incindlink1} is proved.


\end{proo}

By Theorem \ref{t.genindlink1} and Theorem \ref{t.incindlink1}, we  can easily prove Theorem \ref{t.main0}.

In the end of the paper, we provide some remarks about the
 topological equivalence of   NMS flows admitting a given indexed link.
If an NMS flow on a closed $3$-manifold does not  admit any saddle closed orbit, then by FRH decompositions, the closed orbits of this flow consist of an attracting closed orbit and a repelling closed orbit. Here, we refer to this type of flow as  \emph{North-South NMS flow}.

Suppose that a closed $3$-manifold $M$ admits a North-South NMS flow. It is easy to observe that the FRH decomposition of a North-South NMS flow on $M$ is a Heegaard splitting of $M$ such that each handlebody is a solid torus, and the indexed link of this flow consists of the cores of these  handlebodies. Therefore,
$M$ is homeomorphic to $S^3$, $S^1\times S^2$ or a lens spase. 
By Lemma A.1 of Yu \cite{Yu}, it is not difficult to prove that
the North-South NMS flows on a fixed manifold $M$  are topologically equivalent.

Generally, the number of NMS flows admitting a given indexed link $l$ may not be unique, even may be infinite up to topological equivalence.
For example, if $l$ is  a three component unlinked, unknotted link and only one component is index-$k$ ($k=0,1,2$), then there are $8$ NMS flows on $S^3$ with indexed link $l$  up to topological equivalence (by Yu \cite[Proposition 7]{Yu}). In particular, if the number of  index-$1$ knots in $l$ is more than one, then the number of NMS flows admitting the indexed link $l$ may  be infinite up to topological equivalence because of the complexity of heteroclinic trajectories.

%
%
%
%

\vskip 1cm
\noindent  Fangfang Chen

\noindent{\small School of Mathematical Sciences}

\noindent{\small Tongji University, Shanghai 200092, CHINA}

\noindent{\footnotesize{E-mail: fangfangchen\_97@163.com}}
\vskip 2mm

\noindent Bin Yu

\noindent {\small School of Mathematical Sciences}

\noindent{\small Tongji University, Shanghai 200092, CHINA}

\noindent{\footnotesize{E-mail: binyu1980@gmail.com }}

\end{document}